\documentclass[11pt,a4paper]{article}
\usepackage[utf8]{inputenc}
\usepackage{amsmath}
\usepackage{amsfonts}
\usepackage{color}
\usepackage{amssymb}
\usepackage{mathrsfs}
\usepackage{esint}
\usepackage[colorinlistoftodos]{todonotes}
\usepackage[colorlinks=true]{hyperref}
\hypersetup{
	colorlinks=blue,%
	citecolor=red,%
	filecolor=black,%
	linkcolor=blue,%
	urlcolor=blue
}
\usepackage[margin=1.0in]{geometry}

\usepackage[amsmath,thmmarks,hyperref]{ntheorem}
{
	\theoremstyle{nonumberplain}
	\theoremheaderfont{\bfseries}
	\theorembodyfont{\normalfont}
	\theoremsymbol{\mbox{$\Box$}}
	\newtheorem{pf}{Proof.}
}

\numberwithin{equation}{section}
\def\R{\mathbb{R}}
\def\S{\mathbb{S}}
\def\T{\mathbb{T}}
\def\N{\mathbb{N}}

\def\e{\epsilon}
\newcommand{\ud}{\mathrm{d}}

\newtheorem*{conjecture A}{Conjecture A}
\newtheorem{thm}{Theorem}[section]

\newtheorem{lem}{Lemma}[section]
\newtheorem{rem}{Remark}[section]

\newtheorem{step}{Step}
\newtheorem{thm A}{Theorem A}
\newtheorem*{convention}{\textbf{Convention}}
\newtheorem{cor}{Corollary}[section]
\usepackage{upgreek}
\usepackage{graphicx}

\makeatletter

\newdimen\bibspace
\renewenvironment{thebibliography}[1]{%
	\section*{\refname %or \bibname if you use ``book'' as the documentclass
		\@mkboth{\MakeUppercase\refname}{\MakeUppercase\refname}}%
	\list{\@biblabel{\@arabic\c@enumiv}}%
	{\settowidth\labelwidth{\@biblabel{#1}}%
		\leftmargin\labelwidth
		\advance\leftmargin\labelsep
		\itemsep\bibspace
		\parsep\z@skip     %
		\@openbib@code
		\usecounter{enumiv}%
		\let\p@enumiv\@empty
		\renewcommand\theenumiv{\@arabic\c@enumiv}}%
	\sloppy\clubpenalty4000\widowpenalty4000%
	\sfcode`\.\@m}
{\def\@noitemerr
	{\@latex@warning{Empty `thebibliography' environment}}%
	\endlist}

\makeatother

\makeatletter
\allowdisplaybreaks

\allowdisplaybreaks

\begin{document}
\title{article}

	\title{\bf The moving plane method and the  uniqueness of high order elliptic equation with GJMS operator  } 
	\date{\today}
	\author{\medskip Shihong Zhang 
	}
	
	\renewcommand{\thefootnote}{\fnsymbol{footnote}}
\footnotetext[1]{S. Zhang: dg21210019@smail.nju.edu.cn}
	\maketitle

	%MS+++++++++++++++++++++ Abstract +++++++++++++++++++++++++
	\begin{abstract}
In this paper, we study the following high order elliptic equation involving the GJMS operator:  
		\begin{align*}
		\alpha P_{\mathbb{S}^n}v_{\alpha}+2Q_{g_{\mathbb{S}^n}}=2Q_{g_{\mathbb{S}^n}}e^{nv_{\alpha}}.
		\end{align*}
	We establish that if
		 $\alpha>1$  and $n\geq3$,  or if $\alpha\in (1-\epsilon_0, 1)$ with $n=2m\geq4$, 	then $v_{\alpha}\equiv0$. As an application, we present a new proof of the classical Beckner inequality.
		
	\end{abstract}

		\medskip 
		
		{{\bf $\mathbf{2020}$ MSC:}Primary: 35A23, 35J60, Secondary: 35B33, 35B38}
		

		\section{Introduction}

		    The high order conformally invariant operators are a fundamental topic in conformal geometry. For a closed manifold $(M,g)$, an operator $A_{g}$ is said to be conformally invariant if, under a conformal change of metric $\hat{g}=e^{2w}g$, it satisfies  
		\begin{align*}
			A_{\hat{g}}(\phi)=e^{-bw}A_{g}(e^{aw}\phi),
		\end{align*} 
		where $\phi\in C^{\infty}(M)$ and $a$, $b$ are constants.
		 If $n=2$, the Laplace-Beltrami operator $\Delta_g$  is a basic example of a conformally invariant operator, satisfying  
		\begin{align*}
			\Delta_{\hat{g}}(\phi)=e^{-2w}\Delta_{g}\phi.
		\end{align*}
	For a new conformal metric $ \hat{g}$ on $M$, its Gauss curvature satisfies the following critical exponent elliptic equation: 
		\begin{align*}
			-\Delta_{g}w+K_{g}=K_{\hat{g}}e^{2w}.
		\end{align*}
	For $n \geq 3$, there exists a higher-order conformally invariant operator $P_g^n$, known as the GJMS operator, which transforms under a conformal change of metric as 
		\begin{align*}
		P_{g}^n(w)+2Q_{g}^n=2Q_{\hat{g}}^ne^{nw}.
		\end{align*}
	While the explicit formula for the GJMS operator on a general closed manifold is highly intricate, on the standard sphere $\mathbb{S}^n$, it admits a more explicit representation
		\begin{align}
			P_{\S^n}=\begin{cases}
			\displaystyle 	\prod_{k=0}^{\frac{n-2}{2}}\left(-\Delta_{g_{\S^n}}+k(n-k-1)\right)\qquad\qquad\qquad\qquad\qquad\qquad\,\,\quad\mathrm{for}\,\,\,n\,\,\,\mathrm{even},\\
			\displaystyle \left(-\Delta_{g_{\S^n}}+\left(\frac{n-1}{2}\right)^2\right)^{1/2}\prod_{k=0}^{\frac{n-3}{2}}\left(-\Delta_{g_{\S^n}}+k(n-k-1)\right)\qquad\mathrm{for}\,\,\,n\,\,\,\mathrm{odd}.
			\end{cases}
		\end{align}
A well-known problem related to $ P_{\mathbb{S}^n} $ is the prescribed $ Q $-curvature problem: given a function $ Q $, can we find a solution to the PDE  
\begin{align}\label{Prescribed Q}
	P_{\S^n}w+2Q_{g_{\S^n}}=2Qe^{nw}~?
\end{align}	
In fact, a necessary condition for the solvability of this equation, known as the Kazdan-Warner condition, is given by 
\begin{align}\label{K-W condi}
	\int_{\S^n}\langle\nabla x^{l}, \nabla Q\rangle e^{nw}\ud V_{\S^n}=0, \qquad\mathrm{for}\qquad l=1,2,\cdots,n+1.
\end{align}
In the context of the prescribed $Q$-curvature problem, an intriguing example was presented in \cite{S.Zhang 2}, where the traditional non-degeneracy condition is not satisfied. Nevertheless, the equation \eqref{Prescribed Q} still admits at least one solution, highlighting a scenario where the usual assumptions can be partially relaxed.

To introduce the main topic of this paper, we begin by considering the following geometric functional,
		\begin{align*}
		J_{\alpha}(v)=\frac{n\alpha}{2(n-1)!}\fint_{\S^n} vP_{\S^n}v\ud V_{\S^n}+n\fint_{\S^n} v\ud V_{\S^n}-\log\fint_{\S^n} e^{nv}\ud V_{\S^n},
		\end{align*}
		when $\alpha\geq1$, Onofri inequality \cite{Onofri} $(n=2)$ and Beckner inequality \cite{Beckner} $(n\geq 3)$ states that $J_{\alpha}(v)\geq0$ for any $v\in H^{\frac{n}{2}}(\S^n)$,  and $J_{1}(v)=0$ if and only if $e^{2v}g_{\S^n}=\varphi^{*}g_{\S^n}$ for some conformal transformation $\varphi$ on $\S^n$. To proceed with the discussion, we define
		\begin{align*}
		\mathcal{M}=\left\{v\Big|\int_{\S^n}e^{nv}x^{l}\ud V_{\S^n}=0\quad\mathrm{for}\quad l=1,2,\cdots,n+1\right\}.
		\end{align*}
An intriguing problem is to consider the minimum and the critical points of 
 $J_{\alpha}(\cdot)$ within $\mathcal{M}$. Specifically,
		\begin{align}\label{variation problem}
			J_{\alpha}: =\inf_{v\in\mathcal{M}}J_{\alpha}(v).
		\end{align}
		Using the Kazdan-Warner condition \eqref{K-W condi}, one can prove that all the possible critical points of $J_{\alpha}$ satisfy
		\begin{align}\label{sphere main equ}
		\alpha P_{\S^n}v_{\alpha}+2Q_{g_{\S^n}}=2Q_{g_{\S^n}}e^{nv_{\alpha}},
		\end{align}
		where $Q_{g_{\S^n}}=\frac{(n-1)!}{2}$. Regarding \eqref{variation problem}, we conjecture that the following result holds: 
		\begin{align*}
			J_{\alpha}= 0\qquad\mathrm{for}\qquad \alpha\in[1/2, 1)\,\,?
		\end{align*} 
	Motivated by results established by Aubin \cite{Aubin} and Chang–Yang \cite{Chang and Yang anns}, we propose a stronger conjecture:
		\begin{conjecture A}\label{Conjecture}
			For $n\geq2$ and for any $\alpha\in [1/2, 1)$ the following equation
				\begin{align*}
			\alpha P_{\S^n}v_{\alpha}+2Q_{g_{\S^n}}=2Q_{g_{\S^n}}e^{nv_{\alpha}}
			\end{align*}
			has only trivial solution.
		\end{conjecture A} 
		Regarding the above \eqref{variation problem}, \eqref{sphere main equ}, and Conjecture A, several results are listed below.
		\begin{enumerate}
			\item[(1)] $n=2$. Provided that $\alpha>1$,  then $J_{\alpha}= 0$ and the solutions of \eqref{sphere main equ} are only trivial , you can refer \cite{Chen Li,Chanillo Kiessling}. If $\alpha>\frac{1}{2}$, Aubin \cite{Aubin} who proved that the variational problem \eqref{variation problem} is bounded below and is attained by minimizer $v_{\alpha}$ in $\mathcal{M}$. For $\alpha\in (1-\e_0,1)$, where $\e_0$ is sufficiently small,  Chang and  Yang \cite{Chang and Yang acta} showed that $J_{\alpha}=0$ and the minimizer $v_{\alpha}\equiv 0$. With the axially symmetric assumption, Feldman, Froese, Ghoussoub, and Gui \cite{Feldman J. Froese  R Ghoussoub Gui} have proven that  under $\alpha>\frac{16}{25}-\e_0$, all solutions of \eqref{sphere main equ} are zero;  for $\alpha\geq\frac{1}{2}$, Gui and Wei \cite{Gui Wei}, and Lin \cite{Lin} independently proved that \eqref{sphere main equ} has only trivial solution. For general case, Ghoussoub and Lin \cite{Ghoussoub Lin} showed that if $\alpha>\frac{2}{3}-\e_0$, then the Conjecture A is true. Finally, Gui and Moradifam \cite{Gui Moradifam} completed the proof of  Conjecture A by the sphere covering inequality.
			\item[(2)] $n=2m\geq 4$. Suppose $\alpha>\frac{1}{2}$, Chang and Yang \cite{Chang and Yang anns} showed that variational problem \eqref{variation problem} is bounded below and is achieved by minimizer $v_{\alpha}$. Afterwards, Wei and Xu \cite{Wei and Xu} clarified that $J_{\alpha}=0$ and the minimizer $v_{\alpha}\equiv 0$ for $\alpha\in (1-\e_0,1)$. Recently, for axially symmetric case,  Gui, Hu and Xie \cite{Gui Hu Xie} proved that all solutions of \eqref{sphere main equ} are trivial if $\alpha=\frac{1}{n+1}$ or $\alpha\in (\alpha_n, 1)$ with $n=6, 8$. 
		\end{enumerate}

		Similar to $n=2$, if $\alpha>1$, applying the moving plane method as \cite{S.Zhang}, we prove that $J_{\alpha}=0$ and the solution of \eqref{sphere main equ} are trivial for dimension greater than two.
		\begin{thm}\label{main thm 1}
			Suppose $\alpha>1$ and $n\geq3$, if $v_{\alpha}$ is a solution of 
			\begin{align*}
			\alpha P_{\S^n}v_{\alpha}+2Q_{g_{\S^n}}=2Q_{g_{\S^n}}e^{nv_{\alpha}},
			\end{align*}
			then $v_{\alpha}\equiv0$.
		\end{thm}
	 
	 For $\alpha < 1$, the situation becomes more challenging due to the loss of a uniform maximum principle (i.e., Lemma \ref{maximum principle}) for all $\alpha \in [1/2, 1)$. To overcome this difficulty, we adopt the ``compactness implies uniqueness" technique from \cite{Chang and Yang acta}, which relies on the eigenvalue analysis of the GJMS operator. Inspired by \cite{Malchiodi}, we establish the following compactness theorem for $\alpha > \frac{1}{2} + \varepsilon$, where $\varepsilon > 0$ is sufficiently small.

			\begin{thm}\label{main thm 2}
			Suppose $n=2m\geq 4$, then for any $\e>0$ the following set
			\begin{align*}
			\mathcal{S}_{\e}=\left\{v_{\alpha}\Big|\alpha P_{\S^n}v_{\alpha}+2Q_{g_{\S^n}}=2Q_{g_{\S^n}}e^{nv_{\alpha}}\quad\mathrm{and}\quad \alpha\in \left(\frac{1}{2}+\e,1\right)\right\}
			\end{align*}
		is compact in $C^{\infty}(\S^n)$ topology.
		\end{thm}
	
	As a direct consequence, we can prove Conjecture A for 
	 $\alpha\in (1-\e_0, 1)$, thereby extending the results of Wei and Xu \cite{Wei and Xu}. It is also worth noting that while Conjecture A has been fully resolved in two dimensions by Gui and Moradifam \cite{Gui Moradifam}, their approach heavily relies on counting the number of nodal sets and Bol’s inequality.
	\begin{thm}\label{main thm 2.1}
		Suppose $n=2m\geq 4$, there exist $\e_0>0$ such that for all $\alpha\in (1-\e_0,1)$, the equation
		\begin{align*}
		\alpha P_{\S^n}v_{\alpha}+2Q_{g_{\S^n}}=2Q_{g_{\S^n}}e^{nv_{\alpha}}
		\end{align*}
		has only trivial solution.
	\end{thm}

Lastly, we point out a very interesting corollary, i.e, the value of $\alpha$ which has only trivial solution is open in $(1/2,1)$.
	\begin{thm}\label{main thm 2.2}
		Suppose $n=2m\geq4$, denote
		\begin{align*}
		\Gamma	=\left\{\alpha\in \left(\frac{1}{2},1\right)\bigg|\left\{v_{\alpha}|\alpha P_{\S^n}v_{\alpha}+2Q_{g_{\S^n}}=2Q_{g_{\S^n}}e^{nv_{\alpha}}\right\}=\{0\}\right\}.
		\end{align*}
	 then $\Gamma$ is open in $\left(\frac{1}{2},1\right)$ and $(1-\e_0, 1) \subset \Gamma$
	\end{thm}

	The organization of this paper is as follows. In Section \ref{sec2}, we establish essential estimates for $\left(-\Delta\right)^{m-1} u(x)$. In Subsection \ref{sec3.1}, we prove the radial symmetry of solutions with respect to their critical points in three steps: initiating the moving plane method, showing that $w_{\lambda}(x)$ vanishes in the limiting state, and  the plane must stop the origin. Subsection \ref{sec3.2} focuses on determining the possible number of critical points of the solution. We then establish the uniqueness of solutions on the sphere for $\alpha > 1$.  In Subsection \ref{Sec 3.3}, we present a new proof of Beckner’s inequality using Theorem \ref{main thm 1}. In Subsection \ref{sec 4.1}, we develop a concentration-compactness lemma along with key a priori estimates. Subsection \ref{sec 4.2} classifies higher-order conformally invariant equations under certain integral decay conditions. In Subsection \ref{sec 4.3}, we analyze the possible number of concentration points and prove Theorems \ref{main thm 2}, \ref{main thm 2.1}, and \ref{main thm 2.2}. Additionally, we examine the critical case when $\alpha_i \to \frac{1}{2}$, proving Theorem \ref{critical blow up thm} and illustrating the blow-up behavior of the sequence $\{v_i\}_{i=1}^{+\infty}$.

		\section{Preliminaries}\label{sec2}
	In Sections \ref{sec2} and \ref{sec3}, we establish the theorems for a fixed $\alpha$ and any given solution $v_{\alpha}$, meaning the constants in these sections may depend on $\alpha$. Since we only require the assumption $\alpha > 1$, we omit the subscript $\alpha$ for simplicity.  
	
	First, for any fixed point $p \in \mathbb{S}^n$, we introduce the stereographic projection $I = I_p$ with $p$ as the north pole. More specifically, the mapping $I = I_p: (\mathbb{R}^n, x) \to \mathbb{S}^n \setminus \{p\}$ is given by  
			\begin{align*}
				I_p(x)=\left(\frac{2x_1}{1+|x|^2},\cdots, \frac{2x_n}{1+|x|^2}, \frac{|x|^2-1}{|x|^2+1}\right).
			\end{align*}
			For simplicity, we often omit the subscript $p$, through the conformal change of GJMS operator $P_{\S^n}$ we know
		\begin{align*}
		\alpha\left(\frac{2}{1+|x|^2}\right)^{-n}\left(-\Delta\right)^{\frac{n}{2}}v\circ I(x)+2Q_{g_{\S^n}}=2Q_{g_{\S^n}}e^{nv\circ I(x)}.
		\end{align*}
		Denote $u(x)=v\circ I(x)+\frac{1}{\alpha}\log\frac{2}{1+|x|^2}$, it satisfies
		\begin{align*}
		\alpha\left(\frac{2}{1+|x|^2}\right)^{-n}\left(-\Delta\right)^{\frac{n}{2}}&\left(u(x)-\frac{1}{\alpha}\log\frac{2}{1+|x|^2}\right)+2Q_{g_{\S^n}}\\
		=&2Q_{g_{\S^n}}e^{nu(x)}\left(\frac{2}{1+|x|^2}\right)^{-\frac{n}{\alpha}}.
		\end{align*}
		Then,
		\begin{align}\label{main equation 0}
		\left(-\Delta\right)^{\frac{n}{2}}u(x)=2Q_g(x)e^{nu(x)}\qquad\qquad\mathrm{in}\qquad\R^n,
		\end{align}
		where 
		\begin{align*}
		Q_g(x)=\frac{Q_{g_{\S^n}}}{\alpha}\left(\frac{2}{1+|x|^2}\right)^{n\left(1-1/\alpha\right
			)}.
		\end{align*}
		In this paper, denote $f=Q_ge^{nu}$.
		Recall that $u(x)=v\circ I(x)+\frac{1}{\alpha}\log\frac{2}{1+|x|^2}$, we know
		\begin{align}\label{main condi 1}
			\lim_{|x|\to+\infty}\frac{u(x)}{-\frac{2}{\alpha}\log|x|}=1.
		\end{align}
	Suppose $ I' = I_{-p}: (\mathbb{R}^n, y) \to \mathbb{S}^n \setminus \{-p\}$ is another stereographic projection with $-p$ as the north pole. Then, the smoothness of $v$ around $p$ is equivalent to the smoothness of $v \circ I'(y)$ around $0$. Clearly, we have  
	 $y=(I^{'})^{-1}\circ I(x)=\frac{x}{|x|^2}$ and 
	\begin{align*}
	v\circ I(x)=\left(v\circ I^{'}\right)\circ\left(I^{'}\right)^{-1}\circ I(x)=\left(v\circ I^{'}\right)\left(\frac{x}{|x|^2}\right).
	\end{align*}
	Denote $m=\left[\frac{n+1}{2}\right]\geq 2$ and we let $h(x)=\left(-\Delta\right)^{m-1} u(x)$, we always assume that $p$ is a critical point of $v$.
	
	\begin{lem}\label{expan of u}
		For $|x|\gg1$, if $p$ is a critical point of $v$, we have 
		\begin{align*}
				u(x)=-\frac{2}{\alpha}\log|x|+c+O(|x|^{-2}).
		\end{align*}
	\end{lem}
	\begin{pf}
		Let $\bar{u}(x)=	u(x)+\frac{2}{\alpha}\log|x|=v\circ I(x)+\frac{1}{\alpha}\log\frac{2|x|^2}{1+|x|^2}$, then 
		\begin{align*}
			\lim_{|x|\to+\infty}\bar{u}(x)=v(p)+\frac{1}{\alpha}\log 2:=c.
		\end{align*}
		We consider 
		\begin{align*}
			|x|\left(\bar{u}(x)-c\right)=|x|\left(\left(v\circ I^{'}\right)\left(\frac{x}{|x|^2}\right)-v(p)\right)+|x|\log\frac{|x|^2}{1+|x|^2}.
		\end{align*}
		Noticing that $\left(v\circ I^{'}\right)(0)=v(p)$, $\nabla_y\left(v\circ I^{'}\right)(0)=0$ and $\lim_{|x|\to+\infty}|x|\log\frac{|x|^2}{1+|x|^2}=0$, we have
		\begin{align*}
		\left||x|\left(\bar{u}(x)-c\right)\right|=\left|\left(\nabla_y\left(v\circ I^{'}\right)(0)\cdot\frac{x}{|x|}\right)\right|+O\left(\frac{1}{|x|}\right)=O\left(\frac{1}{|x|}\right).
		\end{align*}
	\end{pf}

The following Lemma \ref{Lem 2.2} is a classical result in the moving plane method, and our proof is partially inspired by Lemma 3.2 in \cite{Chang and Yang}.
		\begin{lem}\label{Lem 2.2}
			Suppose $p$ is a critical point of $v$, then $h(x)$ has the following expansion at infinity:
			\begin{align}\label{asy formula a}
					\begin{cases}
					\displaystyle h(x)=\frac{a_0}{|x|^{2(m-1)}}+\sum_{i, j=1}^{n}\frac{a_{ij}x_ix_i}{|x|^{2m+2}}+O\left(\frac{1}{|x|^{2m+1}}\right),\\
					\displaystyle \partial_ih(x)=-2(m-1)\sum_{i=1}^{n}\frac{a_0x_i}{|x|^{2m}}+O\left(\frac{1}{|x|^{2m+1}}\right),\\
					\displaystyle \partial_i\partial_jh(x)=O\left(\frac{1}{|x|^{2m}}\right),
					\end{cases}
			\end{align}
			where $a_0>0$ and $\left|\nabla^k u(x)\right|=O\left(\frac{1}{|x|^{k}}\right)$.
		\end{lem}
		\begin{pf}
			Since  \begin{align*}
				\partial _i\log\frac{2}{1+|x|^2}&=\frac{-2x_i}{1+|x|^2}\\
				\partial_j\partial _i\log\frac{2}{1+|x|^2}&=-2\left(\frac{\delta_{ij}}{1+|x|^2}-\frac{2x_ix_j}{(1+|x|^2)^2}\right),
			\end{align*}
			hence
			\begin{align*}
				\left(-\Delta\right)\log\frac{2}{1+|x|^2}=\frac{n+(n-2)|x|^2}{(1+|x|^2)^2}=\frac{(n-2)}{|x|^2}+O\left(\frac{1}{|x|^4}\right).
			\end{align*}
			Note that the $\log\frac{2}{1+r^2}$ is analytic with respect to the variable $r=|x|$, when $|x|\gg1$ and $$\left(-\Delta\right)^{m-2}\left(\frac{1}{|x|^2}\right)=\frac{c_0}{|x|^{2(m-1)}}$$ where $c_0=\prod_{i=1}^{m-2}(n-2i-2)\prod_{i=1}^{m-2}(2i)>0$, then \begin{align*}
					\left(-\Delta\right)^{m-1}\log\frac{2}{1+|x|^2}=\frac{(n-2)c_0}{|x|^{2(m-1)}}+O\left(\frac{1}{|x|^{2m}}\right).
			\end{align*}
			Now we continue to calculate $(-\Delta)\left(v\circ I(x)\right)$, 
			\begin{align*}
				\partial_i\left(v\circ I(x)\right)&=\partial_{y_{\alpha}}\left(v\circ I^{'}\right)(y)\frac{\partial y_{\alpha}}{\partial x_i}\\
				&=\partial_{y_{\alpha}}\left(v\circ I^{'}\right)(y)\left(\frac{\delta_{\alpha i}}{|x|^2}-\frac{2x_\alpha x_i}{|x|^4}\right)
			\end{align*}
			and
			\begin{align*}
				\partial_j\partial_i\left(v\circ I(x)\right)&=\partial_{y_{\beta}}\partial_{y_{\alpha}}\left(v\circ I^{'}\right)(y)\left(\frac{\delta_{\alpha i}}{|x|^2}-\frac{2x_\alpha x_i}{|x|^4}\right)\left(\frac{\delta_{\beta j}}{|x|^2}-\frac{2x_\beta x_j}{|x|^4}\right)\\
				&+\partial_{y_{\alpha}}\left(v\circ I^{'}\right)(y)\left(-\frac{2\delta_{\alpha i}x_j}{|x|^4}-\frac{2\delta_{\alpha j} x_i+2x_{\alpha}\delta_{ij}}{|x|^4}+\frac{8x_\alpha x_ix_j}{|x|^6}\right).
			\end{align*}
			Recall that $y=\frac{x}{|x|^2}$ , $\partial_{y_{\alpha}}\left(v\circ I^{'}\right)(0)=0$ and $\partial_{y_{\beta}}\partial_{y_{\alpha}}\left(v\circ I^{'}\right)(y)$ are smooth near $y=0$. Thus, $\left|\nabla^2\left(v\circ I(x)\right)\right|=O\left(\frac{1}{|x|^4}\right)$. Similarly, you can get $\left|\nabla^k\left(v\circ I(x)\right)\right|=O\left(\frac{1}{|x|^{k+2}}\right)$. Finally, we prove the formula
			\begin{align*}
				h(x)=\frac{a_0}{|x|^{2(m-1)}}+\sum_{i, j=1}^{n}\frac{a_{ij}x_ix_i}{|x|^{2m+2}}+O\left(\frac{1}{|x|^{2m+1}}\right).
			\end{align*}
		\end{pf}
		\begin{lem}\label{lem 2.2}
			Suppose $h(x)=\left(-\Delta\right)^{m-1} u(x)$, then 
			\begin{align}\label{lem 2.2.0}
				h(x)=C(n)\int_{\R^n}\frac{Q_g(y)e^{nu(y)}}{|x-y|^{2(m-1)}}dy,
			\end{align}
			for some constant $C(n)$.
		\end{lem}
	\begin{pf}
		Let $w(x)=C(n)\int_{\R^n}\frac{Q_g(y)e^{nu(y)}}{|x-y|^{2(m-1)}}dy$, we claim that
		\begin{align}\label{lem 2.2.00}
			\lim_{|x|\to+\infty}|x|^{2(m-1)}w(x)=C(n)\int_{\R^n}Q_g(y)e^{nu(y)}dy.
		\end{align} 
		Denote $f=Q_ge^{nu}$,  $I=\int_{\R^n}\frac{f(y)}{|x-y|^{2(m-1)}}dy$, for $|x|\gg R$, then
		\begin{align}\label{lem 2.2.1}
		I\geq \int_{|y|<R}\frac{f(y)}{|x-y|^{2(m-1)}}dy\geq \frac{\int_{|y|<R}f(y)dy}{(|x|+R)^{2(m-1)}}.
		\end{align}
		Multiplying $|x|^{2(m-1)}$ and taking $|x|\to\infty$ for \eqref{lem 2.2.1}, we obtain
		\begin{align*}
		\liminf_{|x|\to+\infty}|x|^{2(m-1)}w(x)\geq C(n)\int_{|y|<R}f(y)dy,
		\end{align*}
		then let $R\to\infty$ we get
		\begin{align*}
		\liminf_{|x|\to+\infty}	|x|^{2(m-1)}w(x)\geq C(n)\int_{\R^n}f(y)dy.
		\end{align*} 
		Again, we know
		\begin{align}\label{lem 2.2.2}
		I\leq& \int_{|y|<R}+\int_{|y-x|<|x|/2}+\int_{|y-x|>|x|/2, \,|y|>R}\frac{f(y)}{|x-y|^{2(m-1)}}dy\nonumber\\
		:=&I_1+I_2++I_3.
		\end{align}
		Similarly, there holds
		\begin{align}\label{lem 2.2.3}
		I_1\leq \frac{\int_{|y|<R}f(y)dy}{(|x|-R)^{2(m-1)}}.
		\end{align}
		Noticing that
		\begin{align}\label{asy e}
		Q_g(x)e^{nu(x)}&=\frac{Q_{g_{\S^n}}}{\alpha}\left(\frac{2}{1+|x|^2}\right)^{n\left(1-1/\alpha\right
			)}e^{nv\circ I(x)}\left(\frac{2}{1+|x|^2}\right)^{n/\alpha}\nonumber\\
		&=\frac{Q_{g_{\S^n}}}{\alpha}\left(\frac{2}{1+|x|^2}\right)^{n}e^{nv\circ I(x)}\sim O\left(\frac{1}{|x|^{2n}}\right).
		\end{align} 
		For $I_2$, by \eqref{asy e} we arrive that
		\begin{align}\label{lem 2.2.4}
		I_2\leq \frac{C}{|x|^{2n}}\int_{|y-x|<|x|/2}\frac{1}{|x-y|^{2(m-1)}}dy\leq\frac{C}{|x|^{n+2(m-1)}}.
		\end{align}
		Clearly, 
		\begin{align}\label{lem 2.2.5}
		I_3\leq \frac{C}{|x|^{2(m-1)}}\int_{|y|>R}f(y)dy=\frac{o_{R}(1)}{|x|^{2(m-1)}}.
		\end{align}
		From \eqref{lem 2.2.2}, \eqref{lem 2.2.3}, \eqref{lem 2.2.4} and \eqref{lem 2.2.5}, we conclude that 
		\begin{align*}
		\limsup_{|x|\to+\infty}	|x|^{2(m-1)}w(x)\leq C(n)\int_{|y|<R}f(y)dy+o_{R}(1).
		\end{align*}
	Let $ R \to \infty $, completing the proof of the claim \eqref{lem 2.2.00}. Define$ z(x) = h(x) - w(x) $. Since $h(x)$ admits a well-behaved expansion (see Lemma \ref{Lem 2.2}) and the limit $\lim_{|x|\to+\infty} |x|^{2(m-1)} w(x)$ exists, it follows that 
		\begin{align*}
			z(x)=h(x)-w(x)\in L_{1}=\left\{u\in L^1_{\mathrm{loc}}\bigg|\int_{\R^n}\frac{|u(x)|}{1+|x|^{n+1}}\ud x<+\infty\right\}.
		\end{align*}
	 If $n=2m-1$,	hence we can choose $C(n)$ such that $\left(-\Delta\right)^{\frac{1}{2}}w=2f$, , which leads to
	 \begin{align*}
	 (-\Delta)^{\frac{1}{2}}z=	(-\Delta)^{\frac{1}{2}}h-\left(-\Delta\right)^{\frac{1}{2}}w=	(-\Delta)^{\frac{1}{2}}\circ (-\Delta)^{m-1}u-2f=0.
	 \end{align*}
	  Here, two cases are discussed separately.
		\begin{enumerate}
			\item[(1)] 	If $n=2m$, , hence
			then 
			\begin{align*}
			-\Delta z(x)=0\qquad\mathrm{in}\qquad\R^n\qquad\mathrm{and}\qquad\lim_{|x|\to+\infty} h(x)=\lim_{|x|\to+\infty} w(x)=0.
			\end{align*}
			By the Liouville theorem of harmonic functions, we know $z(x)\equiv C$, but $\lim_{|x|\to+\infty} z(x)=0$, then $z(x)\equiv 0$.
			\item[(2)] If $n=2m-1$, then 
			\begin{align*}
			\left(-\Delta\right)^{\frac{1}{2}} z(x)=0\qquad\mathrm{in}\qquad\R^n\qquad\mathrm{and}\qquad\lim_{|x|\to+\infty} h(x)=\lim_{|x|\to+\infty} w(x)=0.
			\end{align*}
			Since $z(x)$ is a bounded function, set $\tilde{z}(x)=z(x)+\sup_{x\in\R^n}|z(x)|\geq0$ and 
			\begin{align*}
				\left(-\Delta\right)^{\frac{1}{2}} \tilde{z}(x)=0\qquad\mathrm{in}\qquad\R^n.
			\end{align*}
			Through the elliptic estimate of fractional Laplacian, we know
			\begin{align}\label{lem 2.2.9}
				||\tilde{z}||_{C^{1,\beta}(B_1(0))}\leq C(n)\inf_{B_{2}(0)}\tilde{z},
			\end{align}
			for some $\beta\in (0,1)$. The above estimate \eqref{lem 2.2.9}, you can refer  Theorem 2.11 in \cite{Li Xiong}. The standard rescaling technique will imply 
			\begin{align*}
				R||\nabla \tilde{z}||_{L^{\infty}(B_R(0))}\leq C(n)\inf_{B_{2R}(0)}\tilde{z},
			\end{align*}
			let $R\to\infty$, we obtain $\tilde{z}\equiv C$. Thus, $z\equiv C$, combining with $\lim_{|x|\to+\infty} z(x)=0$, so $z\equiv 0$.
		\end{enumerate}

	\end{pf}

		\section{Moving plane method}\label{sec3}
		\subsection{radially symmetric aboutits critical point}\label{sec3.1}
		For $\lambda\in\R$, we define 
		\begin{align*}
			T_{\lambda}&=\{x\in\R^n|x_1=\lambda\},\\
			\Sigma_{\lambda}&=\{x\in\R^n|x_1<\lambda\}.
		\end{align*}
		Let $x^{\lambda}=(2\lambda-x_1, x_2, \cdots, x_n)$ be the reflection point of $x$ with respect to $T_{\lambda}$ and 
		\begin{align*}
			w_{\lambda}(x)=u(x^{\lambda})-u(x).
		\end{align*}
		By the above analysis, for $|x|\gg1$ we know
		\begin{align}\label{asy u 1}
				u(x)=-\frac{2}{\alpha}\log|x|+c+O(|x|^{-2})
		\end{align}
		and 
		\begin{align}\label{asy laplace 1}
			\begin{cases}
			\displaystyle h(x)=\frac{a_0}{|x|^{2(m-1)}}+\sum_{i, j=1}^{n}\frac{a_{ij}x_ix_i}{|x|^{2m+2}}+O\left(\frac{1}{|x|^{2m+1}}\right),\\
			\displaystyle \partial_ih(x)=-2(m-1)\sum_{i=1}^{n}\frac{a_0x_i}{|x|^{2m}}+O\left(\frac{1}{|x|^{2m+1}}\right),\\
			\displaystyle \partial_i\partial_jh(x)=O\left(\frac{1}{|x|^{2m}}\right),
			\end{cases}
		\end{align}
		Then we get
		\begin{align}\label{asy infty}
			\liminf_{x\in\Sigma_{\lambda}, \vert x\vert\to +\infty}w_{\lambda}(x)=\liminf_{x\in\Sigma_{\lambda}, \vert x\vert\to +\infty}-\frac{2}{\alpha}\log\frac{|x^{\lambda}|}{|x|}+O(|x|^{-2})+O(|x^{\lambda}|^{-2})=0
		\end{align}
		and
		\begin{align}\label{asy infty 1}
			\lim_{x\in\Sigma_{\lambda}, \vert x\vert\to +\infty}\left|\nabla^kw_{\lambda}(x)\right|=\lim_{x\in\Sigma_{\lambda}, \vert x\vert\to +\infty}\left|\nabla^ku_{\lambda}(x)-\nabla^ku(x)\right|=0.
		\end{align}
		We move the plane by three steps. Step 1, we start the moving plane at sufficiently negative. Step 2, we prove that $w_{\lambda}\equiv 0$ at limit state. Step 3, we obtain the radial symmetry of the solution above its critical point. Our method is somewhat inspired by \cite{Chang and Yang} and \cite{S.Zhang}.
	\begin{step}
		start to move the plane
	\end{step}
		\begin{lem}\label{move plane lem 1}
			There exist $R_0>0$ and $\lambda^1_{0}<0$ such that 
			\begin{align}\label{start palne 1}
				 h(x)<h(x^{\lambda})\qquad \mathrm{for}\qquad\lambda<\lambda^1_{0}\quad\mathrm{and}\quad|x^\lambda|>R_0.
			\end{align}
			
		\end{lem}
	\begin{pf}
		Since 
		\begin{align*}
			h(x^{\lambda})-h(x)=&a_0\left(\frac{1}{|x^{\lambda}|^{2(m-1)}}-\frac{1}{|x|^{2(m-1)}}\right)+O\left(\frac{1}{|x|^{2m}}\right)+O\left(\frac{1}{|x^{\lambda}|^{2m}}\right)\\
			=:&I+O\left(\frac{1}{|x|^{2m}}\right)+O\left(\frac{1}{|x^{\lambda}|^{2m}}\right),
		\end{align*}
		then we estimate the above term in two cases.
		\begin{enumerate}
			\item[(1)] If $|x|>2|x^{\lambda}|$, then $\frac{1}{|x|}<\frac{1}{2|x^{\lambda}|}$. We have\begin{align*}
			I&\geq\left(1-\frac{1}{2^{2(m-1)}}\right)a_0\frac{1}{|x^{\lambda}|^{2(m-1)}}.
			\end{align*}
			So, if $|x^\lambda|>R_0$, $h(x^{\lambda})>h(x)$.
		
	\item[(2)] If $|x^{\lambda}|<|x|<2|x^{\lambda}|$, then $\frac{1}{|x^{\lambda}|}>\frac{1}{|x|}>\frac{1}{2|x^{\lambda}|}$. Then we know
	\begin{align}\label{st pla 1-1}
		I&\geq \frac{a_0}{|x|^{2m-3}}\left(\frac{1}{|x^{\lambda}|}-\frac{1}{|x|}\right)\geq\frac{c_0(|x|-|x^\lambda|)}{|x|^{2m-1}},
	\end{align}
	clearly,
	\begin{align}\label{st pla 1-5}
	|x|-|x^\lambda|=\frac{|x|^2-|x^{\lambda}|^2}{|x|+|x^\lambda|}\geq \frac{-c\lambda(\lambda-x_1)}{|x|}.
	\end{align}

	From the above estimate \eqref{st pla 1-1} and \eqref{st pla 1-5} there holds
	\begin{align}\label{st pla 1-4}
		h(x^{\lambda})-h(x)&\geq \frac{c_0(|x|-|x^\lambda|)}{|x|^{2m-1}}-\frac{C}{|x|^{2m}}\nonumber\\
		&\geq\frac{-C\lambda(\lambda-x_1)}{|x|^{2m}}-\frac{C}{|x|^{2m}}.
	\end{align}
	
	If $\lambda-x_1\geq\frac{C}{-\lambda}$, then $	h(x^{\lambda})>h(x)$.  	If $\lambda+\frac{C}{\lambda}<x_1<\lambda$, we know $-\lambda\gg1$, the second equality in \eqref{asy laplace 1} will also imply $	h(x^{\lambda})>h(x)$.
\end{enumerate}
	\end{pf}

		\begin{lem}\label{start moving plane lem}
			There $\lambda_{0}:=\lambda_{0}(R_0)<\lambda^1_{0}<0$ such that 
			\begin{align}\label{start palne 2}
			 h(x)<h(x^{\lambda})\qquad \mathrm{for}\qquad\lambda<\lambda_{0}\quad\mathrm{and}\quad x\in\Sigma_{\lambda}.
			\end{align}
		\end{lem}
		\begin{pf}
			By Lemma \ref{move plane lem 1}, we know
			\begin{align*}
			h(x)<h(x^{\lambda})
			\end{align*}
			for $\lambda<\lambda_{0}^{1}$ and $|x^{\lambda}|>R_0$. 	From \eqref{lem 2.2.0} and \eqref{lem 2.2.00} in Lemma \ref{lem 2.2}, we know that if $\lambda<\lambda_{0}<\lambda_{0}^{1}$, and since $|x|>|\lambda|$, then we have
			\begin{align*}
			h(x)=O\left(\frac{1}{|x|^{2(m-1)}}\right)<C(R_0):=\min_{x\in B_{R_0}} h(x)\qquad\mathrm{for}\quad |x^{\lambda}|\leq R_0.
			\end{align*} 
		
				\end{pf}
			Denote 
			\begin{align*}
				\bar{\lambda}=\sup\{\lambda\leq 0|\left(-\Delta \right)^{m-1}w_{\mu}>0, \forall~x\in \Sigma_{\mu} , \mu\leq\lambda\},
			\end{align*}
			from Lemma \ref{start moving plane lem}, we know $\bar{\lambda}\geq \lambda_{0}$.
			
			\begin{step}
				If $\bar{\lambda}<0$, then $w_{\bar{\lambda}}(x)\equiv 0$ for $x\in\Sigma_{\bar{\lambda}}$.
			\end{step}
		
	Throughout Step 2, we assume that $\bar{\lambda}<0$. The following maximum principle lemma is the key ingredient for the moving plane approach.
		
		\begin{lem}[Maximum Principle]\label{maximum principle}
			For $\alpha>1$ , if $w_{\bar{\lambda}}\not\equiv 0$ in $\Sigma_{\bar{\lambda}}$, then
			\begin{align*}
				\left(-\Delta\right)^kw_{\bar{\lambda}}>0\qquad\mathrm{for}\qquad k=0,1,\cdots,m-1
			\end{align*}
			 and $\left(-\Delta \right)^{\frac{n}{2}}w_{\bar{\lambda}}>0$ for any $x\in \Sigma_{\bar{\lambda}}$.
		\end{lem}
	\begin{pf}
		 By the definition of $\bar{\lambda}$,  we know $\left(-\Delta\right)^{m-1}w_{\bar{\lambda}}\geq0$. From \eqref{asy infty 1}, we get
		\begin{align}\label{lem 3.3 boun 1}
			\lim_{|x|\to+\infty}\left(-\Delta\right)^{m-2}w_{\bar{\lambda}}=0\qquad\mathrm{and}\qquad \left(-\Delta\right)^{m-2}w_{\bar{\lambda}}=0\qquad\mathrm{on} \qquad T_{\bar{\lambda}}.
		\end{align} 
		Then, maximum principle implies that $\left(-\Delta\right)^{m-2}w_{\bar{\lambda}}\geq0$ in $\Sigma_{\bar{\lambda}}$. Similarly, you can obtain 
		\begin{align}\label{lem 3.3 boun 1.1}
			\left(-\Delta\right)^kw_{\bar{\lambda}}\geq 0\qquad\mathrm{in}\qquad\Sigma_{\bar{\lambda}}, \qquad\mathrm{for}\qquad k=1,\cdots,m-3.
		\end{align}
		Now we have $-\Delta w_{\bar{\lambda}}\geq0$ and \eqref{asy infty} tells us that
		\begin{align}\label{lem 3.3 boun 2}
				\lim_{|x|\to+\infty}w_{\bar{\lambda}}=0\qquad\mathrm{and}\qquad w_{\bar{\lambda}}=0\qquad\mathrm{on} \qquad T_{\bar{\lambda}}.
		\end{align}
		By strong maximum principle and $w_{\bar{\lambda}}\not\equiv 0$, we get $w_{\bar{\lambda}}>0$ for                                                                                                                                           $x\in \Sigma_{\bar{\lambda}}$.
		Due to $w_{\bar{\lambda}}(x)=u(x^{\bar{\lambda}})-u(x)>0$, we see that
		\begin{align}\label{lem 3.3 right term}
&\left(-\Delta \right)^{\frac{n}{2}}w_{\bar{\lambda}}(x)\nonumber\\
&=\frac{2Q_{g_{\S^n}}}{\alpha}\left(\left(\frac{2}{1+|x^{\bar{\lambda}}|^2}\right)^{n\left(1-1/\alpha\right
	)}e^{nu(x^{\bar{\lambda}})}-\left(\frac{2}{1+|x|^2}\right)^{n\left(1-1/\alpha\right
	)}e^{nu(x)}\right)\nonumber\\
		&>0,
		\end{align}
		the above  inequality follows by $|x^{\bar{\lambda}}|<|x|$. The Lemma \ref{lem 2.2} yields that 
			\begin{align*}
				\left(-\Delta\right)^{m-1}h(x)&=C(n)\int_{\R^n}\frac{f(y)}{|x-y|^{2(m-1)}}dy=C(n)\int_{\Sigma_{\bar{\lambda}}}+\int_{\R^n\backslash \Sigma_{\bar{\lambda}}}\frac{f(y)}{|x-y|^{2(m-1)}}dy\\
				&=C(n)\int_{\Sigma_{\bar{\lambda}}}\frac{f(y)}{|x-y|^{2(m-1)}}dy+C(n)\int_{ \Sigma_{\bar{\lambda}}}\frac{f(y^{\bar{\lambda}})}{|x-y^{\bar{\lambda}}|^{2(m-1)}}dy\\
				&=C(n)\int_{\Sigma_{\bar{\lambda}}}\frac{f(y)}{|x-y|^{2(m-1)}}+\frac{f(y^{\bar{\lambda}})}{|x^{\bar{\lambda}}-y|^{2(m-1)}}dy,
			\end{align*}
				this last equality follows by $|x^{\bar{\lambda}}-y|=|x-y^{\bar{\lambda}}|$. Since 
				\begin{align*}
				\left(-\Delta\right)^{m-1}w_{\bar{\lambda}}(x)=f(x^{\bar{\lambda}})-f(x):=g_{\bar{\lambda}}(x),
				\end{align*}
				then
				\begin{align}\label{lem 3.3 lap w}
				\left(-\Delta\right)^{m-1}w_{\bar{\lambda}}(x)&=C(n)\int_{ \Sigma_{\bar{\lambda}}}\left(\frac{1}{|x^{\bar{\lambda}}-y|^{2(m-1)}}-\frac{1}{|x-y|^{2(m-1)}}\right)f(y)dy\nonumber\\
				&+C(n)\int_{ \Sigma_{\bar{\lambda}}}\left(\frac{1}{|x-y|^{2(m-1)}}-\frac{1}{|x^{\bar{\lambda}}-y|^{2(m-1)}}\right)f(y^{\bar{\lambda}})dy\nonumber\\
				&=C(n)\int_{\Sigma_{\bar{\lambda}}}\left(\frac{1}{|x-y|^{2(m-1)}}-\frac{1}{|x^{\bar{\lambda}}-y|^{2(m-1)}}\right)g_{\bar{\lambda}}(y)dy,
				\end{align}
			combining with $|x-y|<|x^{\bar{\lambda}}-y|$ for $x, y\in\Sigma_{\bar{\lambda}}$ and $g_{\bar{\lambda}}(x)>0$ in $\Sigma_{\bar{\lambda}}$ (this is a direct result by \eqref{lem 3.3 right term}), we conclude that
		 $	\left(-\Delta\right)^{m-1}w_{\bar{\lambda}}(x)>0$. Applying strong maximum principle, you can get 
		\begin{align*}
		\left(-\Delta\right)^kw_{\bar{\lambda}}> 0\qquad\mathrm{in}\qquad\Sigma_{\bar{\lambda}}, \qquad\mathrm{for}\qquad k=1,\cdots,m-2.
		\end{align*}

	\end{pf}

			\begin{lem}\label{lem 3.4}
				If $w_{\bar{\lambda}}(x)\not\equiv 0$, then there exists $\delta_0^1\ll 1$ and $R_1\gg \lambda_{0}$ such that 
				\begin{align}
					\frac{\partial}{\partial x_1}\left(-\Delta\right)^{m-1} u(x)>0\quad\forall~\vert x\vert>R_{1},~ \bar{\lambda}- \delta_{0}^{1}<x_1<\bar{\lambda}+ \delta_{0}^{1}<0.
				\end{align}
			\end{lem}
			\begin{pf}
			
			Let $b=(\bar{\lambda},0,\cdots,0), \bar{x}= x-b$, then
			$\vert x\vert\sim \vert x^{\bar{\lambda}}\vert\sim\vert \bar{x}\vert$, if$\,\, \vert x\vert>R_1\gg \max\{\vert \lambda_{0}\vert,R_{0}\}$. Then we have the following claim.
			
				\textbf{Claim:} There exists $\e_0>0$ such that 
				\begin{align}\label{derivatve lap lem 1}
					\left(-\Delta\right)^{m-1}w_{\bar{\lambda}}> \frac{\e_0(\bar{\lambda}-x_1)}{|\bar{x}|^{2m}}\qquad  \forall~x\in\varSigma_{\bar{\lambda}}\backslash B_{R_{1}}.
				\end{align}
			 We observe that 
				\begin{align}\label{derivatve lap lem a1}
				&	\frac{1}{|x-y|^{2(m-1)}}-\frac{1}{|x^{\bar{\lambda}}-y|^{2(m-1)}}\nonumber\\
				&\geq \frac{1}{|x^{\bar{\lambda}}-y|^{2m-3}}\left(\frac{1}{|x-y|}-\frac{1}{|x^{\bar{\lambda}}-y|}\right)
				=\frac{|x^{\bar{\lambda}}-y|-|x-y|}{|x^{\bar{\lambda}}-y|^{2m-2}|x-y|}\nonumber\\
				&\geq \frac{C(\bar{\lambda}-x_1)(\bar{\lambda}-y_1)}{|x^{\bar{\lambda}}-y|^{2m}}.
				\end{align}
				From \eqref{lem 3.3 lap w} and \eqref{derivatve lap lem a1}, we get
				\begin{align*}
						\left(-\Delta\right)^{m-1}w_{\bar{\lambda}}(x)&=C(n)\int_{\Sigma_{\bar{\lambda}}}\left(\frac{1}{|x-y|^{2(m-1)}}-\frac{1}{|x^{\bar{\lambda}}-y|^{2(m-1)}}\right)g_{\bar{\lambda}}(y)dy\\
						&\geq C (\bar{\lambda}-x_1)\int_{\Sigma_{\bar{\lambda}}}\frac{(\bar{\lambda}-y_1)}{|x^{\bar{\lambda}}-y|^{2m}}g_{\bar{\lambda}}(y)dy\\
						&\geq C(\bar{\lambda}-x_1)\int_{\bar{\lambda}-y_1\geq 1, |x^{\bar{\lambda}}-y|\leq 2|x| }\frac{(\bar{\lambda}-y_1)}{|x^{\bar{\lambda}}-y|^{2m}}g_{\bar{\lambda}}(y)dy\\
						&\geq \frac{C(\bar{\lambda}-x_1)}{|x|^{2m}}\int_{\bar{\lambda}-y_1\geq 1, |x^{\bar{\lambda}}-y|\leq 2|x| }g_{\bar{\lambda}}(y)dy\\
						&\geq \frac{\e_0(\bar{\lambda}-x_1)}{|\bar{x}|^{2m}}.
				\end{align*}
			The last inequality follows by 
			\begin{align*}
				\int_{\bar{\lambda}-y_1\geq 1, |x^{\bar{\lambda}}-y|\leq 2|x| }g_{\bar{\lambda}}(y)dy\to\int_{\bar{\lambda}-y_1\geq 1 }g_{\bar{\lambda}}(y)dy>0\qquad\mathrm{as}\qquad|x|\to\infty.
			\end{align*}
		
			Now, we complete the proof the claim.	
				The estimate \eqref{derivatve lap lem 1} directly implies that
				\[
				\frac{\left(-\Delta\right)^{m-1} w_{\bar{\lambda}}(x_1,x')-(\left(-\Delta\right)^{m-1} w_{\bar{\lambda}})(\bar{\lambda},x')}{x_1-\bar{\lambda}}\leq\frac{-\e_0}{\vert \bar{x}\vert^{2m}} \qquad \mathrm{for~~} x=(x_1,x') \in \Sigma_{\bar{\lambda}}\setminus B_{R_{1}}.
				\]
				Letting $x_1\rightarrow \bar{\lambda}$, we obtain that for $\vert x\vert>R_{1}$
				\[\frac{\partial\left(-\Delta\right)^{m-1} w_{\bar{\lambda}} }{\partial x_1}(\bar{\lambda},x')=-2\frac{\partial\left(-\Delta\right)^{m-1} u}{\partial x_1}(\bar{\lambda},x')\quad\Longrightarrow\quad \frac{\partial\left(-\Delta\right)^{m-1} u}{\partial x_1}(\bar{\lambda},x')\geq\frac{\e_0}{2\vert \bar{x}\vert^{2m}}.
				\]
				By Taylor expansion \eqref{asy laplace 1}, we further choose $0<\delta<\delta_{0}^1\ll 1$ such that
				\[\frac{\partial\left(-\Delta\right)^{m-1} u}{\partial x_1}(x_1,x')\geq\frac{\e_0}{4\vert \bar{x}\vert^{2m}},\]
				for all $\vert x\vert>R_{1}, \bar{\lambda}- \delta_{0}^{1}<x_1<\bar{\lambda}+ \delta_{0}^{1}<0$.

			\end{pf}

				\begin{lem}\label{lem3.5}
				If $\lambda_0<\bar{\lambda}<0, w_{\bar{\lambda}}(x)\not\equiv 0$, then there exist $0<\delta_{0}^{2}\ll \delta_{0}^{1}\ll 1,\,R_{2}>R_{1}$ such that
				\[
				\left(-\Delta\right)^{m-1} w_{\lambda}>0,\quad\forall~ \vert x\vert>R_{2}, \quad x_1<\lambda,\quad\lambda \in (\bar{\lambda},\bar{\lambda}+\frac{\delta_{0}^{2}}{4}).
				\]
			\end{lem}
			\begin{pf}
					It follows from Lemma  \ref{lem 3.4} that for $\lambda \in (\bar{\lambda},\bar{\lambda}+\frac{\delta_{0}^{2}}{4})$ and  $\forall~ \vert x\vert>R_{1}$, $\lambda- \frac{\delta_{0}^{1}}{2}<x_1<\lambda,\quad$
						\[
					h (x^{\lambda})>h (x)\quad \Longrightarrow	\quad \left(-\Delta\right)^{m-1} w_{\lambda}>0.
					\]
						Here $\delta_{0}^2$ is chosen such that $\lambda- \frac{\delta_{0}^{1}}{2}<\bar{\lambda}-\frac{\delta_{0}^{1}}{4}$, it suffices to consider $\vert x\vert>R_{1},  x_1<\bar{\lambda}-\frac{\delta_{0}^{1}}{4}$.
					Notice that $\vert x\vert\sim \vert x^{\bar{\lambda}}\vert\sim\vert \bar{x}\vert\sim\vert x^{\lambda}\vert$ for all  $\vert x\vert>R_{1}$, then  
					\begin{equation}\label{equ3.9}
					\left| \frac{1}{\vert  x^{\bar{\lambda}} \vert^k}-\frac{1}{\vert  x^{\lambda} \vert^k}\right|<\frac{C(\lambda-\bar{\lambda})(\vert\bar{\lambda}-x_1\vert+\vert \lambda\vert)}{\vert \bar{x}\vert^{k+2}}.
					\end{equation}
					By (\ref{equ3.9}) we can estimate
					\[\left|\left(-\Delta\right)^{m-1} w_{\lambda}-\left(-\Delta\right)^{m-1} w_{\bar{\lambda}}\right|\leq
					I+II+III+O\left(\frac{1}{\vert \bar{x}\vert^{2m+1}}\right),\]
					where
					\[
					\begin{aligned}
					I&\leq a_0\left|\frac{1}{\vert  x^{\bar{\lambda}} \vert^{2m-2}}-\frac{1}{\vert  x^{\lambda} \vert^{2m-2}}\right|,\\
					II&\leq\left|\frac{a_{11}(2\lambda-x_1)^2}{\vert x^{\lambda}\vert^{2m+2}} -\frac{a_{11}(2\bar{\lambda}-x_1)^2}{\vert x^{\bar{\lambda}}\vert^{2m+2}}\right|,\\
					III&\leq2\sum_{i=2}^{n}\left| \frac{a_{1i}(2\lambda-x_1) x_i}{\vert x^{\lambda}\vert^{2m+2}}-\frac{a_{1i}(2\bar{\lambda}-x_1) x_i}{\vert x^{\bar{\lambda}}\vert^{2m+2}}\right|,\\
					IV&\leq \sum_{i,j=2}^{3}\left|a_{ij}x_ix_j\right|\left|\frac{1}{\vert x^{\bar{\lambda}}\vert^{2m+2}}-\frac{1}{\vert  x^{\lambda} \vert^{2m+2}}\right|.
					\end{aligned}
					\]
					For $I$, by \eqref{equ3.9}, we know 
					\begin{align}\label{Lem 3.5 a1}
						|I|\leq \frac{C(\lambda-\bar{\lambda})(\vert\bar{\lambda}-x_1\vert+\vert \lambda\vert)}{\vert \bar{x}\vert^{2m}}.
					\end{align}
					The second term follows by 
					\begin{align}\label{Lem 3.5 a2}
						|II|\leq& \left|\frac{a_{11}(2\lambda-x_1)^2}{\vert x^{\lambda}\vert^{2m+2}} -\frac{a_{11}(2\bar{\lambda}-x_1)^2}{\vert x^{\lambda}\vert^{2m+2}}\right|+\left|\frac{a_{11}(2\bar{\lambda}-x_1)^2}{\vert x^{\lambda}\vert^{2m+2}} -\frac{a_{11}(2\bar{\lambda}-x_1)^2}{\vert x^{\bar{\lambda}}\vert^{2m+2}}\right|\nonumber\\
						\leq&\frac{C|\lambda-\bar{\lambda}|\left(\vert x^{\lambda}\vert+\vert x^{\bar{\lambda}}\vert\right)}{|\bar{x}|^{2m+2}}+\frac{C|x^{\bar{\lambda}}|^2(\lambda-\bar{\lambda})(\vert\bar{\lambda}-x_1\vert+\vert \lambda\vert)}{|\bar{x}|^{2m+4}}\nonumber\\
						\leq& \frac{C|\lambda-\bar{\lambda}|}{|\bar{x}|^{2m+1}}.
					\end{align}
					Similarly, you can obtain
					\begin{align}\label{Lem 3.5 a3}
						|III|\leq& 2\sum_{i=2}^{n}\left| \frac{a_{1i}(2\lambda-x_1) x_i}{\vert x^{\lambda}\vert^{2m+2}}-\frac{a_{1i}(2\bar{\lambda}-x_1) x_i}{\vert x^{\lambda}\vert^{2m+2}}\right|+\left| \frac{a_{1i}(2\bar{\lambda}-x_1) x_i}{\vert x^{\lambda}\vert^{2m+2}}-\frac{a_{1i}(2\bar{\lambda}-x_1) x_i}{\vert x^{\bar{\lambda}}\vert^{2m+2}}\right|\nonumber\\
						\leq& \frac{C|\lambda-\bar{\lambda}||\bar{x}|}{|\bar{x}|^{2m+2}}+\frac{C|x^{\bar{\lambda}}||x|(\lambda-\bar{\lambda})(\vert\bar{\lambda}-x_1\vert+\vert \lambda\vert)}{|\bar{x}|^{2m+4}}\nonumber\\
						\leq& \frac{C|\lambda-\bar{\lambda}|}{|\bar{x}|^{2m+1}}.
					\end{align}
					Finally, from \eqref{equ3.9} you have
					\begin{align}\label{Lem 3.5 a4}
					|IV|\leq\frac{C|x|^2(\lambda-\bar{\lambda})(\vert\bar{\lambda}-x_1\vert+\vert \lambda\vert)}{|\bar{x}|^{2m+4}}	\leq \frac{C|\lambda-\bar{\lambda}|}{|\bar{x}|^{2m+1}}.
					\end{align}
					The \eqref{Lem 3.5 a1}, \eqref{Lem 3.5 a2}, \eqref{Lem 3.5 a3} and \eqref{Lem 3.5 a4} yield that
					\[
					\left|\left(-\Delta\right)^{m-1} w_{\lambda}-\left(-\Delta\right)^{m-1} w_{\bar{\lambda}}\right|<\frac{C(\lambda-\bar{\lambda})(\vert\bar{\lambda}-x_1\vert+\vert \lambda\vert+C)}{\vert \bar{x}\vert^{2m}}+O\left(\frac{1}{\vert \bar{x}\vert^{2m+1}}\right).
					\]
					Thus, we have
					\[
					\begin{aligned}
				\left(-\Delta\right)^{m-1} w_{\lambda}&\geq \left(-\Delta\right)^{m-1} w_{\bar{\lambda}}-\left|\left(-\Delta\right)^{m-1}w_{\lambda}-\left(-\Delta\right)^{m-1} w_{\bar{\lambda}}\right|\\
					&\geq \frac{\e_0 (\bar{\lambda}-x_1)}{\vert \bar{x}\vert^{2m}}-\frac{C(\lambda-\bar{\lambda})(|\bar{\lambda}-x_1|+\vert \lambda\vert+C)}{\vert \bar{x}\vert^{2m}}+O\left(\frac{1}{\vert \bar{x}\vert^{2m+1}}\right).
					\end{aligned}
					\]
					Since  $\bar{\lambda}-x_1> \frac{\delta_{0}^{1}}{4}$,  you can choose $\delta_{0}^{2}\ll \delta_{0}^{1},\,\,\vert \bar{x}\vert>R_{2}\gg R_{1}$,\, such that $	\left(-\Delta\right)^{m-1} w_{\lambda}>0$.

			\end{pf}

			\begin{lem}\label{lem3.6}
				With the same assumption as in Lemma \ref{lem3.5}, there exists $0<\delta_{0}^{3}\ll\delta_{0}^{2}$ such that $	\left(-\Delta\right)^{m-1}w_{\lambda}>0,\,\,\forall~ x\in\varSigma_{\lambda},\,\, \lambda\in(\bar{\lambda},\bar{\lambda}+\delta_{0}^{3}).$
			\end{lem}
			\begin{pf}
				Suppose not, by Lemma \ref{lem3.5} we can choose 
				$$\lambda_{i}\in (\bar{\lambda},\bar{\lambda}+\frac{\delta_{0}^{2}}{4})\rightarrow \bar{\lambda},\quad x_{i}\in B_{R_{2}}(0)\cap\varSigma_{\lambda_{i}}$$
				such that
				$$	\left(-\Delta\right)^{m-1} w_{\lambda_{i}}(x_{i})=
				\min\limits_{\varSigma_{\lambda_{i}}}	\left(-\Delta\right)^{m-1} w_{\lambda_{i}}\leq0.$$
				Up to a subsequence, we can assume
				\[\lim\limits_{i\to \infty}x_{i}=x_{0},\quad\nabla 	\left(-\Delta\right)^{m-1} w_{\bar{\lambda}}(x_{0})=0,\quad	\left(-\Delta\right)^{m-1} w_{\bar{\lambda}}(x_{0})\leq0.\]
				On the other hand, we can apply the Lemma\,\ref{lem 3.4} to show
				\begin{align}\label{lem 3.6 lim 1}
				\left(-\Delta\right)^{m-1} w_{\bar{\lambda}}>0,\quad x\in \varSigma_{\bar{\lambda}};\quad\left(-\Delta\right)^{\frac{n}{2}} w_{\bar{\lambda}}>0,\quad x\in \varSigma_{\bar{\lambda}}.
				\end{align}
				The \eqref{lem 3.6 lim 1} forces $x_{0}\in T_{\bar{\lambda}}$, there are two cases will happen.
				\begin{enumerate}
					\item[(1)] If $n=2m$,  the classical Hopf lemma gives
					\[
					\frac{\partial\left(-\Delta\right)^{m-1}w_{\bar{\lambda}} }{\partial x_1}(x_{0})<0.
					\]
					This is a contradiction with $\nabla 	\left(-\Delta\right)^{m-1} w_{\bar{\lambda}}(x_{0})=0$.
					\item[(2)] If $n=2m-1$, then \eqref{lem 3.3 lap w} will imply
					\begin{align*}
						\frac{\partial\left(-\Delta\right)^{m-1}w_{\bar{\lambda}} }{\partial x_1}(x_{0})=-4C(n)(m-1)\int_{\Sigma_{\bar{\lambda}}}\frac{(\bar{\lambda}-y_1)g_{\bar{\lambda}}(y)}{|x_0-y|^{2m}}dy<0.
					\end{align*}
				\end{enumerate}
				Thus, we complete the argument for all dimension.
			\end{pf}

			Therefore, from Lemmas \ref{lem3.6}, we finish the proof of \textbf{Step 2}.

				\begin{step}
			$\bar{\lambda} = 0$, and $u$ is radially symmetric about the origin.
			\end{step}
		
			Otherwise, if $\bar{\lambda}<0$, then \textbf{Step 2} implies
			$$u_{\bar{\lambda}}(x)=u(x).$$
			$\\\
		$
			This together with
			\[
			\begin{aligned}
			\left(-\Delta\right)^{\frac{n}{2}}u(x)&=\frac{2Q_{g_{\S^n}}}{\alpha}\left(\frac{2}{1+|x|^2}\right)^{n\left(1-1/\alpha\right
				)}e^{nu(x)}\\
			\left(-\Delta\right)^{\frac{n}{2}}u_{\bar{\lambda}}(x)&=\frac{2Q_{g_{\S^n}}}{\alpha}
			\left(\frac{2}{1+|x^{\bar{\lambda}}|^2}\right)^{n\left(1-1/\alpha\right
				)}e^{nu_{\bar{\lambda}}(x)}
			\end{aligned}
			\]
			implies 
			\[(1+|x|^2)^{n\left(1/\alpha-1\right
				)}=
		(1+|x^{\bar{\lambda}}|^2)^{n\left(1/\alpha-1\right
			)}\Longrightarrow \quad \bar{\lambda}=0.\]
			This is a contradiction. By the definition of $\bar{\lambda}$, we know $\left(-\Delta\right)^{m-1}w_{\bar{\lambda}}\geq0$, by the similar analysis as \eqref{lem 3.3 boun 1} and \eqref{lem 3.3 boun 2}, we know $w_{\bar{\lambda}}\geq0$. Combining with $\bar{\lambda}=0$, we know 
			\begin{align*}
				u(-x_1,x_2,\cdots,x_n)\geq u(x_1,x_2,\cdots,x_n)\quad\mathrm{for} \quad x\in\Sigma_{0}=\{x\in\R^n|x_1<0\}.
			\end{align*}
			Similarly, you can also get
			\begin{align*}
			u(-x_1,x_2,\cdots,x_n)\leq u(x_1,x_2,\cdots,x_n)\quad\mathrm{for} \quad x\in\Sigma_{0}=\{x\in\R^n|x_1<0\}.
			\end{align*}
			
			Recall the notation $u(x) = v \circ I(x) + \frac{1}{\alpha} \log \frac{2}{1 + |x|^2}$ and $I = I_p : (\mathbb{R}^n, x) \to \mathbb{S}^n \setminus \{p\}$, where $p$ is a critical point of $v$. We know that $u$ being radially symmetric about the origin is equivalent to $v$ being radially (or axially) symmetric about $p$.

				\begin{thm}\label{radial symmetric critical point thm}
				Suppose $\alpha>1$ and $n\geq3$, if $v$ is the solution of 
				\begin{align*}
				\alpha P_{\S^n}v+2Q_{g_{\S^n}}=2Q_{g_{\S^n}}e^{nv},
				\end{align*}
				then $v$ is radially symmetric about its critical point.
			\end{thm}
		\begin{pf}
			The theorem follows by the The Step 1, Step 2 and Step 3.
		\end{pf}
		\subsection{uniqueness result for $\alpha>1$}\label{sec3.2}
				\begin{lem}\label{lem4.1}
				If $v\not\equiv 0$, then there exist only two antipodal critical points of $v$.
			\end{lem}
			\begin{pf}
				First, we note that $v\not\equiv 0$ is equivalent to $v\not\equiv C$. We assert that there exist at least one minimal point $p$ and one maximum point $q$ of $v$ such that  $d_{\S^n}(p,q)=\pi$. Without loss of generality, we assume  $p=S$.
				
				If $r=d_{\S^n}(p,q)\leq\frac{\pi}{2}$, then
			$z\in \partial B_{r}(p)\cap \partial B_{r}(q)\not =\emptyset$. By Theorem \ref{radial symmetric critical point thm}, $v$ is radially symmetric with respect to $d_{\S^n}(\cdot,p)$ and $d_{\S^n}(\cdot,q)$, then
			\[
			v(z)=	v(p)=	v(q).
			\]
			This contradicts the assumption that $v\not\equiv C$.
			
			If $r=d_{\S^n}(p,q)>\frac{\pi}{2}$ and $q\not =N$, then $r'=\pi-r=d_{\S^n}(N,q)<\frac{\pi}{2}$. In $B_{r^{'}}(N)$, by Theorem \ref{radial symmetric critical point thm}, we know $v$ take constant (its maximum) on $\partial B_{r^{'}}(N)$ , thus two cases may occur.
			
			(1) If there exists $x_{0}\in B_{r'}(N)$ such that 
			\[
			v(x_0)=\min\limits_{B_{r'}(N)}v<v(q),
			\]
			then  $x_{0}$  is also a critical point and $d_{\S^n}(x_{0},q)<\frac{\pi}{2}$. A similar argument yields $v(q)=v(x_{0})$. A contradiction!
			
			(2)\,\,$v\equiv v(q)$ in $B_{r'}(N)$. We define  
			\[
			r_0 = \sup \{ r \mid v \equiv v(N) = v(q) \text{ in } B_r(N) \}.
			\]  
			Clearly, from the previous discussion, we have $ r' \leq r_0 \leq \pi/2 $.  
			
			Now, consider any point $ \tilde{q} \in \partial B_{r_0}(N) $. Since $ \tilde{q} $ is also a maximum point of $ v$, and $ v$ is radially symmetric with respect to $ d_{\mathbb{S}^n}(\cdot, \tilde{q}) $, it follows that for any $ r \leq r_0 $, the intersection $ \partial B_r(\tilde{q}) \cap B_{r_0}(N) \neq \emptyset $. This implies that $ v \equiv v(N) $ in $ B_{r_0}(\tilde{q}) $.  
			
			However, since $ \tilde{q}$ was chosen arbitrarily on $\partial B_{r_0}(N)$, this leads to a contradiction: we can find a larger radius (for example) $ r_0' = \frac{3}{2} r_0$  such that $ v \equiv v(N) = v(q)$ in $ B_{r_0'}(N)$, contradicting the definition of $ r_0 $.

				Next we claim that there is no other critical point except for the above two critical points of $v$. Otherwise, there exists a third critical point $x_{0} $ of $v$, then either $d_{\S^n}(x_{0},q)\leq \frac{\pi}{2}$ or $d_{\S^n}(x_{0},p)\leq \frac{\pi}{2}$ holds. A similar argument also shows that $v\equiv 0$. This contradicts the assumption.
			\end{pf}
			
				\begin{thm}
				Suppose $\alpha>1$ and $n\geq3$, if $v$ is a solution of 
				\begin{align*}
				\alpha P_{\S^n}v+2Q_{g_{\S^n}}=2Q_{g_{\S^n}}e^{nv},
				\end{align*}
				then $v\equiv 0$.
			\end{thm}
			\begin{pf}
				By contradiction if $v\not\equiv 0$, then it follows from Lemma \ref{lem4.1} and Theorem \ref{radial symmetric critical point thm} that $v$ is increasing  along the great circle from the minimum point (say, $S$) to the maximum point (say, $N$). This yields that
				\[\nabla v\cdot\nabla x^{n+1}\geq 0,
				\]
				here $\nabla:=\nabla_{g_{\S^n}}$.
				From the  Kazdan-Warner condition (see Lemma \ref{Kazdan-Warner lem}), we have
				\[
				\int_{\S^n}\langle\nabla v, \nabla x^{n+1}\rangle_{g_{\S^n}} \ud V_{\S^n}=n\int_{\S^n} x^{n+1}v\ud V_{\S^n}=0.
				\]
				These facts together imply that $\nabla_{g_{\S^n}}v=0$, and thus $v\equiv 0$.
			\end{pf}
		
		\subsection{A new proof of Beckner inequality}\label{Sec 3.3}
	As an interesting application, we present a new proof of the Beckner inequality, following the strategy originally developed by Osgood, Phillips, and Sarnak \cite{OPS}. Additionally, we note that the classical Beckner inequality can also be derived from the fractional case, as shown by Xiong \cite{Xiong}.
	
	First, we introduce an almost sharp Beckner-type inequality, which is a direct consequence of the concentration compactness phenomenon and the Moser-Trudinger type inequality. For further details, refer to \cite{Branson&Chang&Yang, Hang, Hang1}.
		
		\begin{lem}\label{Almost sharp lem}
			For any $\e>0$, there exists $C_{\e}$ such that if $v\in H^{\frac{n}{2}}(\S^n)$, we have
			\begin{align*}
				\log \fint_{\S^n}e^{nv}   \ud V_{{\S^n}}\leq \frac{n(1+\e)}{2(n-1)!}\fint_{\S^n} vP_{\S^n}v\ud V_{\S^n}+n\fint_{\S^n} v\ud V_{\S^n}+C_{\e}.
			\end{align*}
		\end{lem}
		
		Now we can prove the Beckner inequality by Theorem \ref{main thm 1} and Lemma \ref{Almost sharp lem}.
		
		\begin{thm}
			For any $v\in H^{\frac{n}{2}}(\S^n)$, we have \begin{align*}
					\log \fint_{\S^n}e^{nv}   \ud V_{{\S^n}}\leq \frac{n}{2(n-1)!}\fint_{\S^n} vP_{\S^n}v\ud V_{\S^n}+n\fint_{\S^n} v\ud V_{\S^n}
			\end{align*}
			and the equality holds if and only if  $e^{2v}g_{\S^n}=\varphi^{*}g_{\S^n}$ (up to a constant) for some conformal transformation $\varphi$ on $\S^n$.
		\end{thm}
		\begin{pf}
			For fixed $\e>0$, by Lemma \ref{Almost sharp lem}, we have $J_{1+\e}\geq -C_{\e}$, thus we can choose a minimizing sequence $\{v_i\}_{i=1}^{+\infty}$ such that $\fint_{\S^n}  v_i\ud V_{{\S^n}}=0$ and 
			\begin{align*}
				J_{1+\e}(v_i)\to J_{1+\e}:=\inf_{v\in H^{\frac{n}{2}}(\S^n)}J_{1+\e}(v).
			\end{align*}
		On the one hand, we can see that if $\fint_{\S^n}  v \ud V_{{\S^n}}=0$, then
			\begin{align}\label{Sec 3.3 formula a}
			 C_1(n)||v||_{H^{\frac{n}{2}}(\S^n)}\leq \fint_{\S^n}   vP_{\S^n}v\ud V_{{\S^n}}\leq C_2(n)||v||_{H^{\frac{n}{2}}(\S^n)},
			\end{align}
			the above inequality follows by the spherical harmonic expansion or you can also see the Page 18 in \cite{Hang1}.
			On the other hand, we also know $J_{1+\e/2}(v_i)\geq -C_{\e/2}$, then
			\begin{align}\label{Sec 3.3 formula b}
			\fint_{\S^n}   v_iP_{\S^n}v_i\ud V_{{\S^n}}\leq C(\e).
			\end{align}
			Now, we get $||v_i||_{H^{\frac{n}{2}}(\S^n)}\leq C(\e)$ from \eqref{Sec 3.3 formula a} and \eqref{Sec 3.3 formula b}, without loss of generality, we can assume $v_i\rightharpoonup v_{\e}\in H^{\frac{n}{2}}(\S^n)$ and  $\fint_{\S^n}   v_{\e}\ud V_{{\S^n}}=0$, then
			\begin{align*}
				J_{1+\e}\leq J_{1+\e}(v_{\e})\leq \liminf_{i\to+\infty} J_{1+\e}(v_i):=J_{1+\e}.
			\end{align*}
			In conclusion, $v_{\e}$ is a minimizer of functional  $J_{1+\e}(\cdot)$, it satisfies
			\begin{align*}
				(1+\e)P_{\S^n}v_{\e}+2Q_{g_{\S^n}}=2Q_{g_{\S^n}}\frac{e^{nv_{\e}}}{\fint_{\S^n}e^{nv_{\e}}   \ud V_{{\S^n}}}.
			\end{align*}
			By Theorem \ref{main thm 1}, we get $v_{\e}\equiv\fint_{\S^n}  v_{\e} \ud V_{{\S^n}}=0$. So far, we have proved that for any $\e>0$, 
			\begin{align*}
				\log \fint_{\S^n}e^{nv}   \ud V_{{\S^n}}\leq \frac{n(1+\e)}{2(n-1)!}\fint_{\S^n} vP_{\S^n}v\ud V_{\S^n}+n\fint_{\S^n} v\ud V_{\S^n},
			\end{align*}
			letting $\e\to 0$, we get the Beckner inequality. Lastly, if $v$ is a minimizer of $J_{1}(\cdot)$, then 
			\begin{align*}
			P_{\S^n}v+2Q_{g_{\S^n}}=2Q_{g_{\S^n}}\frac{e^{nv}}{\fint_{\S^n}e^{nv}   \ud V_{{\S^n}}},
			\end{align*}
			by the classification theorem in \cite{Wei and Xu 1}, we complete the proof.
		\end{pf}

		\section{The compactness of solutions for $\alpha\in\left(\frac{1}{2}+\e,1\right)$}	\label{sec 4}
		In the following sections, we consider $\alpha<1$ and $n=2m$, all notations will add $\alpha$  are subscriptd with $\alpha$ to emphasize $\alpha$ and the relationship between them.
		Recall the main equation,
		\begin{align*}
			 P_{\S^n}v_{\alpha}+\frac{2}{\alpha}Q_{g_{\S^n}}=\frac{2}{\alpha}Q_{g_{\S^n}}e^{nv_{\alpha}}.
		\end{align*}
	Roughly speaking,	if the compactness theorem is not true, then there exists $\alpha_i\to \alpha$ such that $||v_{\alpha_i}||_{C^0(\S^n)}\to+\infty$. Since $\fint_{\S^n}e^{nv_{\alpha_i}}\ud V_{\S^n}=1$ and $\alpha\in (\frac{1}{2}, 1]$, there is at most one blow up point $x_1$ and
		\begin{align}\label{simp blow up}
	\liminf_{i\to+\infty}\frac{e^{nv_{\alpha}}}{|\S^n|}\ud V_{\S^n}\geq\alpha\delta_{p_0}\qquad\mathrm{and}\qquad \lim_{i\to+\infty}\int_{\S^n\backslash B_{r}(x_1)}e^{nv_{\alpha_i}}\ud V_{\S^n}=0
		\end{align}
		However Kazdan-Warner condition yields that
		\begin{align}\label{Kazdan-Warner a}
		\fint_{\S^n}e^{nv_{\alpha}}x^l\ud V_{\S^n}=0\qquad\mathrm{for}\qquad l=1,2,\cdots,n+1.
		\end{align}
		If we choose $x
		_1$ as the north pole, then \eqref{simp blow up} will imply
		\begin{align*}
		\liminf_{i\to+\infty}\fint_{\S^n}e^{nv_{\alpha}}x^{n+1}\ud V_{\S^n}\geq\alpha>0,
		\end{align*}
		this is  a contradiction with \eqref{Kazdan-Warner a}.
		\begin{convention}
			In the following, let $x^{l}$ be the coordinate function on $\S^n$, i.e,
			\begin{align*}
				\S^n=\{x\in\R^{n+1}|\left(x^1\right)^2+\left(x^2\right)^2\cdots\left(x^{n+1}\right)^2=1\}.
			\end{align*}
			And,  $x$, $x_k$ are the points on $\S^n$ or $\R^n$. We also let $c_n=2^{n-1}\left(\frac{n-2}{2}\right)!\pi^{n/2}=Q_{g_{\S^n}}|\S^n|$  throughout the Section 4.
		\end{convention}

		\subsection{Some preliminaries}\label{sec 4.1}
	In this subsection, $x$ and $y$ be the point on $\S^n$ and $|x-y|:=|x-y|_{g_{\S^n}}$ be the distance on $\S^n$.	And, $B_r(x)$ be the geodesic ball on $\S^n$ or general manifold. For closed manifold $(M,g)$, the Green's function of $P_g^{n}$ satisfies
	\begin{align}\label{Green func}
		\left|G_g(x,y)-\frac{1}{c_n}\log\frac{1}{|x-y|}\right|\leq C\qquad\mathrm{for}\qquad x\not=y
	\end{align}
		and 
		\begin{align}\label{derivative Green func}
			\left|\nabla_{x}^kG_{g}(x,y)\right|\leq\frac{C}{|x-y|^k}\qquad\mathrm{for}\qquad k=1,\cdots,n-1.
		\end{align}
Where \eqref{Green func}	follows from \cite{Chang and Yang anns} and \cite{Malchiodi}, for \eqref{derivative Green func} you can refer \cite{Aubin 1}.

		\begin{lem}\label{Kazdan-Warner lem}
			If $\alpha\not=1$ and $v_{\alpha}$ solves 
				\begin{align}\label{Kazdan-Warner lem: equa}
			P_{\S^n}v_{\alpha}+\frac{2}{\alpha}Q_{g_{\S^n}}=\frac{2}{\alpha}Q_{g_{\S^n}}e^{nv_{\alpha}},
			\end{align}then 
			\begin{align}\label{direc assu}
			\int_{\S^n}x^lv_{\alpha}\ud V_{\S^n}=0\quad\mathrm{and}\quad	\int_{\S^n}x^le^{nv_{\alpha}}\ud V_{\S^n}=0.\quad\mathrm{for}\quad l=1,2,\cdots,n+1,
			\end{align}
			where $x^l$ is the $l$-th coordinate function on $\S^n$.
		\end{lem}
		\begin{pf}
			We rewrite \eqref{Kazdan-Warner lem: equa} as
			\begin{align*}
			P_{\S^n}v_{\alpha}+2Q_{g_{\S^n}}=&2Q_{g_{\S^n}}e^{nv_{\alpha}}+\left(2-\frac{2}{\alpha}\right)Q_{g_{\S^n}}\\
			=&2Q_{g_{\S^n}}\left(\left(1-\frac{1}{\alpha}\right)e^{-nv_{\alpha}}+1\right)e^{nv_{\alpha}}.
			\end{align*}
			From Kazdan-Warner condition \eqref{K-W condi} and $\alpha\not=1$, we have
			\begin{align*}
			0=&\int_{\S^n}\nabla x^l\cdot\nabla e^{-nv_{\alpha}}e^{nv_{\alpha}}\ud V_{\S^n}=-n\int_{\S^n}\nabla x^l\cdot\nabla v_{\alpha}\ud V_{g_{\S^4}}\\
			=&n^2\int_{\S^n} x^lv_{\alpha}\ud V_{\S^n}
			\end{align*}
			Multiply $x^l$ and integrate both sides for \eqref{sphere main equ}, then 
			\begin{align*}
			\int_{\S^n}x^le^{nv_{\alpha}}\ud V_{\S^n}=0,
			\end{align*}
			the final result follows by $P_{\S^n}x^l=n!x^l$.
		\end{pf}

		        \begin{lem}\label{priori estimate lem}
		        	For any small $r>0$ and $x_0\in \S^n$, if $n-kp>0$ then
		        	\begin{align}\label{priori estimate a}
		        	\int_{B_r(x_0)}|\nabla^kv_{\alpha}|^p\ud V_{\S^n}\leq Cr^{n-kp}\qquad\mathrm{for}\qquad k=1,\cdots,n-1.
		        	\end{align}
		        	Especially, if $k=1,2,\cdots,n-1$, we have
		        	\begin{align*}
		        	||\nabla^k v_{\alpha}||_{L^p(\S^n)}\leq C\qquad\mathrm{for}\qquad 1\leq p<\frac{n}{k}.
		        	\end{align*}
		        \end{lem}
		        \begin{pf}
		        	Let $f_{\alpha}=\frac{2}{\alpha}Q_{g_{\S^n}}(e^{nv_{\alpha}}-1)$, then $\int_{\S^n}|f_{\alpha}|\leq C$. Denote $\bar{v}_{\alpha}=\fint_{\S^n} v_{\alpha}\ud V_{\S^n}$, then we can see
		        	\begin{align*}
		        	v_{\alpha}(x)=\bar{v}_{\alpha}+\int_{\S^n}G_{\S^n}(x,y)f_{\alpha}(y)\ud V_{\S^n}(y)
		        	\end{align*}
		        	and
		        	\begin{align*}
		        	|\nabla_{x}^kG_{\S^n}(x,y)|\leq\frac{C}{|x-y|^{k}}\qquad\mathrm{for}\qquad k=1,\cdots,n-1.
		        	\end{align*}
		        	We know
		        	\begin{align*}
		        	\int_{B_r(x_0)}|\nabla^kv_{\alpha}|^p(x)\ud V_{\S^n}(x)\leq& C\int_{B_r(x_0)}\left(\int_{\S^n}\frac{||f_{\alpha}||_{L^1(\S^n)}}{|x-y|^{k}}\frac{|f_{\alpha}(y)|}{||f_{\alpha}||_{L^1(\S^n)}}\ud V_{\S^n}(y)\right)^p\ud V_{\S^n}(x)\\
		        	\leq& C\int_{B_r(x_0)}\int_{\S^n}\frac{||f_{\alpha}||^p_{L^1(\S^n)}}{|x-y|^{kp}}\frac{|f_{\alpha}(y)|}{||f_{\alpha}||_{L^1(\S^n)}}\ud V_{\S^n}(y)\ud V_{\S^n}(x)\\
		        	\leq&C\sup_{y\in\S^n}\int_{B_r(x_0)}\frac{1}{|x-y|^{kp}}\ud V_{\S^n}(x)\int_{\S^n}|f_{\alpha}(y)|\ud V_{\S^n}(y)\\
		        	\leq& C\int_{B_{2r}(y)}\frac{1}{|x-y|^{kp}}\ud V_{\S^n}(x)\leq Cr^{n-kp},
		        	\end{align*}
		        	the second inequality follows by the Jensen inequality.
			\end{pf}

		The following concentration compactness theorem has appeared in the Remark 3.2 of \cite{Malchiodi}, and we give the details for the convenience of the reader.

		\begin{lem}\label{concentration compactness lem}
			Suppose $||f_i||_{L^1(M, g_i)}\leq C$ and a sequence of closed manifolds $(M, g_i)$  converges to closed manifold $(M, g)$, if $v_{i}$ satisfies
			 \begin{align*}
				P_{g_i}v_{i}=f_i\qquad\mathrm{in}\qquad
				M.
			\end{align*}
			Then up to a subsequence, either there exists $p>n$ such that
			\begin{align}\label{concentration compactness 1}
				\int_{M}e^{p(v_i-\bar{v}_i)}\ud V_{g_i}\leq C,
			\end{align}
			where $\bar{v}_i=\fint_{M}v_i\ud V_{g_i}$or there exist concentration points $x_1,x_2,\cdots,x_{j_0}$ such that for any $r>0$, we have
			\begin{align}\label{concentration compactness 2}
				\liminf_{i\to+\infty}\int_{B_r(x_k)}|f_i|\ud V_{g_{i}}\geq c_n,
			\end{align}
			where $c_n=2^{n-1}\left(\frac{n-2}{2}\right)!\pi^{n/2}$ and $B_r(x_k)$ is the geodesic ball with respect to $(M,g)$.
		\end{lem}
\begin{rem}
If $x_1,x_2,\cdots,x_{j_0}$  are the concentration points of $\{v_i\}_{i=1}^{+\infty}$, then  $x_1,x_2,\cdots,x_{j_0}$  are also concentration points for its any subsequence. By replacing $\{v_i\}_{i=1}^{+\infty}$ with its subsequence, which we still call it by $\{v_{i}\}_{i=1}^{+\infty}$,  we can always assume that the number of concentration points $ j_0 $ are the maximal (since $||f_i||_{L^1(M, g_i)}\leq C$, the number of concentration point of its any subsequence has uniformly upper bound), i.e , for any subsequence $\{v_i^{'}\}_{i=1}^{+\infty}\subset\{v_i\}_{i=1}^{+\infty}$, the number of concentration point of $\{v_i^{'}\}_{i=1}^{+\infty}$ does not exceed  $j_0$.
\end{rem}	
\begin{pf}

	Suppose the second case does not happen, namely, for any $x\in 
	M$ there exist $r_x$ and $\delta_{x}>0$ such that
	\begin{align*}
		\int_{B_{r_x}(x)}|f_i|\ud V_{g_{i}}\leq c_n-\delta_{x}\qquad\mathrm{for}\qquad i\gg 1,
	\end{align*}
where $B_{r_x}(x)$ is the geodesic ball with respect to $(M,g)$. Then, there exist $(M, g)\subset\cup_{k=1}^{N}B_{r_{x_k}/2}(x_k)$, for simplicity, denote $B_k:=B_{r_{x_k}}(x_k)$ and $B_{k/2}:=B_{r_{x_k}/2}(x_k)$. Hence, for any $x\in B_{k/2}$, we know
\begin{align*}
	v_i(x)-\bar{v}_i=\int_{M}G_{i}(x,y)f_{i}(y)\ud V_{g_{i}}(y),
\end{align*}
where $G_{i}(x,y)$ is the Green's function of $P_{g_i}$ and $G_i\to G$. Since $G(x,y)=\frac{1}{c_n}\log\frac{1}{|x-y|}+O(1)$ for $|x-y|<1$, we get

\begin{align*}
	&\int_{B_{k/2}}e^{p(v_i(x)-\bar{v}_i)}\ud V_{g_{i}}(x)\\
	\leq& \int_{B_{k/2}}\exp\left(p\int_{B_{k}}G_{i}(x,y)\left|f_i(y)\right|\ud V_{g_{i}}(y)\right.
	\\
	&\left.+p\int_{M\backslash B_k}G_{i}(x,y)\left|f_i(y)\right|\ud V_{g_{i}}(y)\right)\ud V_{g_{i}}(x)\\
	\leq&C(r_{x_k})\int_{B_{k/2}}\int_{B_{k}}\exp\left(p G_{i}(x,y)||f_i||_{L^1(B_k, g_i)}\right)\frac{\left|f_i(y)\right|}{||f_i||_{L^1(B_k, g_i)}}\ud V_{g_{i}}(y)\ud V_{g_{i}}(x)\\
	\leq&C(r_{x_k})\int_{B_{k/2}}\int_{B_{k}}\left(\frac{1}{|x-y|}\right)^{\frac{p(c_n-\delta_{x_k})}{c_n}}\frac{\left|f_i(y)\right|}{||f_i||_{L^1(B_k, g_i)}}\ud V_{g_{i}}(y)\ud V_{g_{i}}(x).
\end{align*}
	If we choose $n<p<\frac{nc_n}{c_n-\delta_{x_k}}$ for any $k=1,2,\cdots, N$, then
\begin{align*}
\int_{B_{k/2}}e^{p(u_i(x)-\bar{u}_i)}\ud V_{\S^n}(x)\leq& C(r_{x_k})\int_{B_{2r_{x_k}}(x_k)}\left(\frac{1}{|x-x_k|}\right)^{\frac{p(c_n-\delta_{x_k})}{c_n}}\ud V_{g_{i}}(x)\\
\leq& C(r_{x_k}).
\end{align*}
Although $B_{2r_{x_k}}(x_k)$ is the geodesic ball w.r.t $(M,g)$, the above integral is also uniformly bounded when $i\gg1$.
Finally, we obtain
\begin{align*}
\int_{M}e^{p(u_i(x)-\bar{u}_i)}\ud V_{g_{i}}(x)\leq\sum_{k=1}^{N}\int_{B_{k/2}}e^{p(u_i(x)-\bar{u}_i)}\ud V_{g_{i}}(x)\leq C.
\end{align*}

\end{pf}

\begin{cor}\label{cor 4.1}
	Suppose  $v_{i}:=v_{\alpha_i}$ satisfies
\begin{align*}
P_{\S^n}v_{i}+\frac{2}{\alpha}Q_{g_{\S^n}}=\frac{2}{\alpha}Q_{g_{\S^n}}e^{nv_{i}},
\end{align*}
 and the second case in Lemma \ref{concentration compactness lem} occurs, then for  $1\leq k\leq n$,
 \begin{align}
\lim_{r\to 0} \liminf_{i\to+\infty}\int_{B_r(x_k)}e^{nv_i}\ud V_{\S^n}\geq\frac{\alpha|\S^n|}{2}.
 \end{align}
\end{cor}
\begin{pf}
	Since $f_i=\frac{2}{\alpha_i}Q_{g_{\S^n}}e^{nv_{i}}-\frac{2}{\alpha_i}Q_{g_{\S^n}}$, then for any small $r>0$ we have
	\begin{align*}
		\int_{ B_{r}(x_k)}|f_i|\ud V_{\S^n}=\frac{2}{\alpha_i}Q_{g_{\S^n}}	\int_{ B_{r}(x_k)}e^{nv_{i}}\ud V_{\S^n}+o_{r}(1),
	\end{align*}
	Lemma \ref{concentration compactness lem} will yield that
	\begin{align*}
		\frac{2}{\alpha}Q_{g_{\S^n}}\liminf_{i\to+\infty}	\int_{ B_{r}(x_k)}e^{nv_{i}}\ud V_{\S^n}+o_{r}(1)\geq c_n=Q_{g_{\S^n}}|\S^n|,
	\end{align*}
	letting $r\to 0$, we complete the argument.
\end{pf}

\begin{lem}\label{compactness lem}
	Suppose  $v_{i}:=v_{\alpha_i}$ satisfies
	\begin{align*}
	P_{\S^n}v_{i}+\frac{2}{\alpha}Q_{g_{\S^n}}=\frac{2}{\alpha}Q_{g_{\S^n}}e^{nv_{i}}
	\end{align*}
	and  there exists $p>n$ such that
	 \begin{align}\label{compactness lem a0}
	\int_{\S^n}e^{p(v_i-\bar{v}_i)}\ud V_{\S^n}\leq C.
	\end{align}
	Then there exists $C$ such that
	\begin{align}\label{compactness formula}
	||v_{i}||_{L^{\infty}(\S^n)}\leq C \qquad\mathrm{and}\qquad||\nabla^kv_{i}||_{L^{\infty}(\S^n)}\leq C\qquad\mathrm{for}\qquad k\geq 1.
	\end{align}
\end{lem}
\begin{pf}
 Firstly,
	Holder inequality will imply
	\begin{align}\label{compactness lem a}
	1=e^{n\bar{v}_i}\fint_{\S^n}e^{n(v_i-\bar{v}_i)}\ud V_{\S^n}\leq e^{n\bar{v}_i}\left(\fint_{\S^n}e^{p(v_i-\bar{v}_i)}\ud V_{\S^n}\right)^{\frac{n}{p}}\leq Ce^{n\bar{v}_i}.
	\end{align}
	We also notice that 
	\begin{align}\label{compactness lem b}
	e^{n\bar{v}_i}\leq\fint_{\S^n}e^{nv_i}\ud V_{\S^n}=1,
	\end{align}
	then \eqref{compactness lem a} and \eqref{compactness lem b} will imply
	\begin{align}\label{compactness lem c}
	|\bar{v}_i|=\left|\fint_{\S^n}v_i\ud V_{\S^n}\right|\leq C.
	\end{align}
	Hence, \eqref{compactness lem a0} and \eqref{compactness lem c} will imply that
	\begin{align}\label{compactness lem d}
	\int_{\S^n}e^{pv_i}\ud V_{\S^n}\leq C\qquad\mathrm{for}\quad\mathrm{some}\qquad p>n.
	\end{align}
	Recall $f_i=\frac{2}{\alpha_i}Q_{g_{\S^n}}e^{nv_{i}}-\frac{2}{\alpha_i}Q_{g_{\S^n}}$, then $||f_i||_{L^q(\S^n)}\leq C$ for some $q>1$, this follows by \eqref{compactness lem d}.
	But we also know
	\begin{align}\label{compactness lem e}
	\left|v_i(x)-\bar{v}_i\right|=&\left|\int_{\S^n}G_{\S^n}(x,y)f_{i}(y)\ud V_{\S^n}(y)\right|\nonumber\\
	\leq&\left(\int_{\S^n}\left|G_{\S^n}(x,y)\right|^{q^{'}}\ud V_{\S^n}(y)\right)^{\frac{1}{q^{'}}}\left(\int_{\S^n}|f_{i}(y)|^q\ud V_{\S^n}(y)\right)^{\frac{1}{q}}\\
	\leq& C,\nonumber
	\end{align}
	where $1/q+1/q^{'}=1$.
	The \eqref{compactness lem c} and \eqref{compactness lem e} will imply $||v_i||_{L^{\infty}(\S^n)}\leq C$.  Then $||f_i||_{L^{\infty}(\S^n)}\leq C$, we obtain
	\begin{align*}
	\nabla^kv_i(x)=\int_{\S^n}\nabla_{x}^kG_{\S^n}(x,y)f_{i}(y)\ud V_{\S^n}(y)
	\end{align*}
	and 
	\begin{align*}
	\left|\nabla_{x}^kG_{\S^n}(x,y)\right|\leq\frac{C}{|x-y|^k},
	\end{align*}
	 we conclude that for $1\leq k\leq n-1$,
	\begin{align*}
	||\nabla^kv_{i}||_{L^{\infty}(\S^n)}\leq C.
	\end{align*}
The high order elliptic estimate (see the appendix of \cite{Chen Lu Zhu}) will imply
\begin{align*}
	||\nabla^kv_{i}||_{L^{\infty}(\S^n)}\leq C\qquad\mathrm{for}\qquad k\geq 1.
\end{align*}
\end{pf}

\subsection{A classification result for high order conformally invariant equation}\label{sec 4.2}
In this subsection, we will classify the following equation
	\begin{align}\label{euclidean equation}
\begin{cases}
\displaystyle \left(-\Delta\right)^{\frac{n}{2}}v=(n-1)!e^{nv}\qquad\mathrm{in}\qquad\R^n,\\
\displaystyle
\int_{\R^n}e^{nv(x)}\ud x<+\infty,
\end{cases}
\end{align}
under decay condition: for any $1\leq k\leq n-2$, there exist $\e_k>0$ such that
\begin{align}\label{euclidean equation: decay condition}
	\int_{B_R(0)}|\nabla^k v|\leq CR^{n-\e_k}\qquad\mathrm{for\,\,\,\,\,any}\quad R>0.
\end{align}

For $n=2m\geq4$, we also note that the solution of equation \eqref{euclidean equation} may not be standard bubble without condition \eqref{euclidean equation: decay condition},  which is a big difference between high dimensions and two dimensions. There are many classic results on this issue, see \cite{Chen Li,Lin1,Martinazzi,Wei and Xu 1}. Although there are many works dedicated to proving the integral equation, our primary focus here is on establishing it under the integral constraint \eqref{euclidean equation: decay condition}. For an alternative proof of the integral equation, we refer the reader to \cite{Wei and Ye}.
\begin{lem}\label{classification small lem}
	Suppose $R>0$, there holds
	\begin{align*}
		\int_{B_R(0)}\frac{1}{|x-y|^k}\ud x\leq C(n,k)\left(R^{n-k}\chi_{B_{2R}(0)}(y)+\frac{R^n}{|y|^k}\chi_{B^{c}_{2R}(0)}(y)\right).
	\end{align*}
\end{lem}
\begin{pf}
	If $y\in B_{2R}(0)$, then
	\begin{align}\label{Lemma 4.5 a1}
		\int_{B_R(0)}\frac{1}{|x-y|^k}\ud x\leq\int_{B_{3R}(y)}\frac{1}{|x-y|^k}\ud x\leq C(n)R^{n-k}.
	\end{align}
		For $y\in B^{c}_{2R}(0)$, then $|x-y|\geq |y|-|x|\geq 2R-R=R$. Clearly, 
		\begin{align*}
		|x-y|\leq|x|+|y|\leq R+|y|\leq \frac{3|y|}{2}
		\end{align*}
		and
		\begin{align*}
			|y|\leq |x|+|x-y|\leq R+|x-y|\leq 2|x-y|.
		\end{align*}
		Consequently, we deduce that
		\begin{align}\label{Lemma 4.5 a2}
			\int_{B_R(0)}\frac{1}{|x-y|^k}\ud x\leq \frac{C(n,k)R^n}{|y|^k}.
		\end{align}
		The final result follows by \eqref{Lemma 4.5 a1} and \eqref{Lemma 4.5 a2}.
\end{pf}

In the following, we list two useful results above high order elliptic equations, see \cite{Martinazzi}.

	\begin{lem}[\protect{Pizzetti formula}]\label{Pizzetti}
	Suppose $\Delta^mh=0$ in $B_{2R}(x_0)$, then 
	\begin{align*}
	\fint_{B_{R}(x_0)}h(y)\ud y=\sum_{i=0}^{m-1}c_iR^{2i}\Delta^ih(x_0),
	\end{align*}
	where $c_0=1, c_i=\frac{n}{n+2i}\frac{(n-2)!!}{(2i)!!(2i+n-2)!!}$ for $i\geq1$.
\end{lem}
\begin{lem}\label{high elliptic}
	If  $\Delta^m h=0$  in $B_4(0)$, then for any $\beta\in[0,1)$, $p\geq1$ and $k\geq1$, we have
	\begin{align*}
	||h||_{W^{k,p}(B_1(0))}\leq& C(k,p)||h||_{L^1(B_4(0))},\\
	||h||_{C^{k,\beta}(B_1(0))}\leq& C(k,p)||h||_{L^1(B_4(0))}.
	\end{align*}
\end{lem}	

Inspired by \cite{Lin1,Wei and Xu 1,S.Zhang 1},  we now begin to compare the solution of \eqref{euclidean equation} with the integral
\begin{align*}
	\frac{(n-1)!}{c_n}\int_{\R^n}\log\frac{|y|}{|x-y|}e^{nv(y)}\ud y,
\end{align*}
where $c_n=2^{n-1}\left(\frac{n-2}{2}\right)!\pi^{n/2}$. Finally,  we will prove that the solution of \eqref{euclidean equation} is equal to the above integral plus a constant.
	\begin{lem}\label{upper bound}
Setting 
\begin{align*}
w(x)=\frac{(n-1)!}{c_n}\int_{\R^n}\log\frac{|x-y|}{|y|}e^{nv(y)}\ud y,
\end{align*}then for $|x|\gg1$,
	\begin{align*}
	w(x)\leq a\log|x|+C,
	\end{align*}
	where $a=\frac{(n-1)!}{c_n}\int_{\R^n}e^{nv(y)}\ud y$.
\end{lem}
\begin{pf}
	We decompose $\R^n$ to $\R^n=A_1\cup A_2\cup A_3$, where 
	\begin{align*}
	A_1=&\left\{y||y|<|x|/2\right\},\quad A_2=\left\{y||y-x|<|x|/2\right\},\\
	A_3=&\left\{y||y|>|x|/2,|y-x|>|x|/2\right\}.
	\end{align*}
	In $A_2$, we have $|x-y|<\frac{|x|}{2}<|y|$, then $\log\frac{|x-y|}{|y|}\leq0$. In $A_1$,  we know $\frac{|x|}{2}\leq|x-y|\leq \frac{3|x|}{2}$
	\begin{align}\label{Lemma 4.8 a1}
	\int_{A_1}\log\frac{|x-y|}{|y|}e^{nv(y)}\ud y\leq& \log|x|\int_{A_1}e^{nv(y)}\ud y+C\int_{A_2}f(y)\ud y-\int_{A_1}\log|y|e^{nv(y)}\ud y\nonumber\\
	\leq&\log|x|\int_{A_1}e^{nv(y)}\ud y+C-\max_{B_1}e^{nv}\int_{|y|<1}\log|y|\ud y\nonumber\\
	\leq&\log|x|\int_{A_1}e^{nv(y)}\ud y+C.
	\end{align}
	In $A_3$, we get $|x-y|\leq |y|+|x|\leq3|y|$ and $|y|\leq|x-y|+|x|\leq 3|x-y|$. Then,
	\begin{align}\label{Lemma 4.8 a2}
	\int_{A_3}\log\frac{|x-y|}{|y|}e^{nv(y)}\ud y\leq& C\int_{A_3}e^{nv(y)}\ud y.
	\end{align}
	So, from \eqref{Lemma 4.8 a1} and \eqref{Lemma 4.8 a2} we conclude that
	\begin{align*}
	w(x)\leq a\log|x|+C.
	\end{align*}
\end{pf}

\begin{lem}
	If $v$ is a solution of 
	\begin{align*}
	\begin{cases}
	\displaystyle \left(-\Delta\right)^{\frac{n}{2}}v=(n-1)!e^{nv}\qquad\mathrm{in}\qquad\R^n,\\
	\displaystyle
	\int_{\R^n}e^{nv(x)}\ud x<+\infty,
	\end{cases}
	\end{align*}
	then
	\begin{align*}
		v(x)=\frac{(n-1)!}{c_n}\int_{\R^n}\log\frac{|y|}{|x-y|}e^{nv(y)}\ud y+p(x),
	\end{align*}
	where $p(x)$ is polynomial and $\deg p\leq n-2$.
\end{lem}
\begin{pf}

	Let $h(x)=v(x)+w(x)$, for any fixed $x_0\in\R^n$ by the elliptic estimate of Lemma \ref{high elliptic}  we have
\begin{align}\label{Lem 4.3.a1}
|\nabla^{n-1}h(x_0)|\leq& \frac{C}{R^{n-1}}\fint_{B_{R}(x_0)}|h|\nonumber\\
=&-\frac{C}{R^{n-1}}\fint_{B_{R}(x_0)}h	+\frac{2C}{R^{n-1}}\fint_{B_{R}(x_0)}h^{+}.
\end{align}
Thanks to Lemma \ref{Pizzetti}, we know
\begin{align}\label{Lem 4.3.a2}
\frac{C}{R^{n-1}}\fint_{B_{R}(x_0)}h=O(R^{-1}).
\end{align}
Now we focus on the second part,  applying Lemma \ref{upper bound} and Jensen inequality we obtain
\begin{align}\label{Lem 4.3.1}
\frac{1}{R^{n-1}}\fint_{B_{R}(x_0)}h^{+}\leq& \frac{1}{R^{n-1}}\fint_{B_{R}(x_0)}v^{+}+C\frac{\log R}{R^{n-1}}\nonumber\\
\leq& \frac{C}{R^{n-1}}\log\fint_{B_{R}(x_0)}e^{nv^{+}}+C\frac{\log R}{R^{n-1}}.
\end{align}
Since 
\begin{align}\label{Lem 4.3.2}
\int_{B_{R}(x_0)}e^{nv^{+}}\leq \int_{B_{R}(x_0)}e^{nv}+1
\leq  C+ CR^{n}
\end{align}
for $R\gg1$, from \eqref{Lem 4.3.1} and \eqref{Lem 4.3.2} we arrive 
\begin{align}\label{Lem 4.3.3}
\frac{1}{R^{n-1}}\fint_{B_{R}(x_0)}h^{+}\leq \frac{C}{R^{n-1}}+\frac{C\log R}{R^{n-1}}\to 0,
\end{align}
as $R\to\infty$. Hence, \eqref{Lem 4.3.a1}, \eqref{Lem 4.3.a2} and \eqref{Lem 4.3.3} yield that $\nabla^{n-1}h\equiv 0$, then we complete the argument.
\end{pf}

\begin{thm}\label{classification thm }
	If $v$ is a solution of 
	\begin{align*}
		\begin{cases}
		\displaystyle \left(-\Delta\right)^{\frac{n}{2}}v=(n-1)!e^{nv}\qquad\mathrm{in}\qquad\R^n,\\
		\displaystyle
		\int_{\R^n}e^{nv(x)}\ud  x<+\infty,
		\end{cases}
	\end{align*}
	 and for any $k=1,2,\cdots,n-2$ there exists $\e_k>0$ such that
	 \begin{align}\label{classification thm assumption}
	 	\int_{B_R(0)}|\nabla^k v|\leq CR^{n-\e_k}\qquad\mathrm{for\,\,\,\,\,any}\quad R>0.
	 \end{align}
	 Then,
	  \begin{align*}
	 	v(x)=\log\frac{2\lambda}{\lambda^2+|x-x_0|^2}\qquad\mathrm{and}\qquad \int_{\R^n}e^{nv}=|\S^n|.
	 \end{align*}
\end{thm}
\begin{rem}
If  $n=2$, the constraint in \eqref{classification thm assumption} automatically disappears, making the above theorem a notable generalization of Chen and Li's result in \cite{Chen Li}.
\end{rem}

\begin{pf}
	We will claim :
	\begin{align}\label{classification thm a}
	v(x)=\frac{(n-1)!}{c_n}\int_{\R^n}\log\frac{|y|}{|x-y|}e^{nv(y)}\ud y+C.
	\end{align}
	If the claim is true, in fact, $v(x)=-a(1+o(1))\log|x|$ for $|x|\gg1$, you can refer to the Theorem 2 in \cite{S.Zhang 1}. Then,  the classification of solutions was completed by \cite{Wei and Xu 1}, see also \cite{Martinazzi}. Suppose not, $n-2\geq \deg p=k_0\geq1$.
	Thus, 
	\begin{align*}
	\int_{B_R(0)}	|\nabla^{k_0} v(x)|\geq\int_{B_R(0)}|\nabla^{k_0}p(x)|-\int_{B_R(0)}|\nabla^{k_0}w(x)|, 
	\end{align*}
	here $|\nabla^k_0 v|=\sum_{l_1+\cdots+l_m=k_0}|\partial_{l_1}\cdots\partial_{l_m}v|$.
	Since we know
	\begin{align}\label{classification thm b}
		\int_{B_{R}(0)}|\nabla^{k_0}w(x)|\ud x\leq& C(n)\int_{B_{R}(0)}\int_{\R^n}\frac{e^{nv(y)}}{|x-y|^{k_0}}\ud y\ud x\nonumber\\
		\leq& C(n)\int_{\R^n}\left(\chi_{B_{2R}(0)}(y)R^{n-{k_0}}+\frac{R^n}{|y|^{k_0}}\chi_{B^{c}_{2R}(0)}(y)\right)e^{nv(y)}\ud y\nonumber\\
		\leq& CR^{n-{k_0}},
	\end{align}
	the second inequality follows by Lemma \ref{classification small lem}. And, $\int_{B_{R}(0)}|\nabla^{k_0}p(x)|\ud x\sim  R^n$, then \eqref{classification thm a} and \eqref{classification thm b} will give
	\begin{align*}
\int_{B_R(0)}	|\nabla^{k_0} v(x)|\geq C(n)R^n-C(n,k_0)R^{n-k_0}\geq C(n)R^n\qquad\mathrm{for}\qquad R\gg 1.
	\end{align*}
	But this is a contradiction with the assumption \eqref{classification thm assumption}.
\end{pf}

\subsection{Blow up analysis}\label{sec 4.3}
In this subsection, we will prove that the second case in Lemma \ref{concentration compactness lem} will not occur. If not, we assume that $\{x_k\}_{k=1}^{j_0}$ is the concentration points  of $\{v_{i}\}_{i=1}^{+\infty}$ and $j_0$ is maximal.  We also note that many of the results in this subsection hold in the sense of difference by one subsequence. We will prove:
\begin{enumerate}
	\item[(1)]  From Corollary \ref{cor 4.1}, we know that for a fixed $k\in\{1,2,\cdots,j_0\}$,
	  \begin{align}\label{con vom lower}
	\liminf_{i\to+\infty}\int_{B_r(x_k)}e^{nv_i}\ud V_{\S^n}+o_{r}(1)\geq\frac{\alpha|\S^n|}{2}\qquad\mathrm{for\,\,\,\,\,any}\qquad r>0,
	\end{align}
	then  up to a subsequence,
	\begin{align}\label{concentratiuon volum a}
	\lim_{i\to+\infty}\int_{B_r(x_k)}e^{nv_i}\ud V_{\S^n}\geq\alpha|\S^n|\qquad\mathrm{for\,\,\,\,\,any}\qquad r>0.
	\end{align}
	Now we choose $r<r_0:=\frac{1}{2}\min_{m\not=l}\mathrm{dist}(x_m,x_l)$, noticing that $\alpha\in (1/2, 1]$ and combining with $\int_{\S^n}e^{nv_i}\ud V_{\S^n}=|\S^n|$, then $j_0=1$.
	\item[(2)] Assuming that $x_1$ be the concentration point,  our second argument is that up to a subsequence, for any $r>0$
	\begin{align}\label{sec 4.2:non concentration volum vanish}
	\lim_{i\to+\infty}\int_{\S^n\backslash B_{r}(x_1)}e^{nv_i}\ud V_{\S^n}=0.
	\end{align}
\end{enumerate}
However,  \eqref{concentratiuon volum a} and \eqref{sec 4.2:non concentration volum vanish} can lead to a contradiction with the Kazdan-Warner condition:
\begin{align*}
\int_{\S^n}x^{n+1}e^{nv_{i}}\ud V_{\S^n}=0.
\end{align*}

We also note that this type of blow-up analysis can be significantly simplified using known results. In the context of blow-up analysis on closed manifolds, open domains, and integral equations, there exist numerous seminal works by Hyder, Martinazzi, Druet, Ndiaye, and others. Here, we mention only a few representative references: \cite{Adimurthi&Robert&Struwe,Druet&Robert,Hyder,Hyder&Mancini&Martinazzi,Malchiodi,Martinazzi 1,Martinazzi 2,Ndiaye,Ndiaye 1}.

\begin{lem}
	 Suppose $k\in\{1,2,\cdots,j_0\}$,
	\begin{align*}
	\liminf_{i\to+\infty}\int_{B_r(x_k)}e^{nv_i}\ud V_{\S^n}+o_r(1)\geq\frac{\alpha|\S^n|}{2}\qquad\mathrm{for\,\,\,\,\,any}\qquad r>0.
	\end{align*}
	For any $\rho\in \left(0,\frac{\alpha|\S^n|}{4}\right)$, then there exist $\{r_i\}_{i=1}^{+\infty}$ and $\{x_{k,i}\}_{i=1}^{+\infty}$such that
	\begin{align*}
		\lim_{i\to+\infty}r_i=0\qquad\mathrm{and}\qquad\lim_{i\to+\infty}x_{k,i}=x_k^{'}\in \bar{B}_{r_0}(x_k)
	\end{align*}
	and 
	\begin{align*}
	\int_{B_{r_i}(x_{k,i})}e^{nv_i}\ud V_{\S^n}=\max_{x\in B_{r_0}(x_k)}\int_{B_{r_i}(x)}e^{nv_i}\ud V_{\S^n}=\rho.
	\end{align*}
\end{lem}

\begin{pf}

If \eqref{con vom lower} holds, for fixed $i$, we know that
\begin{align*}
	\max_{x\in B_{r_0}(x_k)}\int_{B_{r}(x)}e^{nv_i}\ud V_{\S^n}
\end{align*}
is continuous above $r$.
Hence, for $\rho\in\left(0,\frac{\alpha|\S^n|}{4}\right)$, then there exist $\{r_i\}_{i=1}^{+\infty}$ such that
\begin{align*}
\int_{B_{r_i}(x_{k,i})}e^{nv_i}\ud V_{\S^n}=\max_{x\in B_{r_0}(x_k)}\int_{B_{r_i}(x)}e^{nv_i}\ud V_{\S^n}=\rho.
\end{align*}
Then we claim that $\lim_{i\to+\infty} r_i=0$. Otherwise, there exists a subsequence $\{r_i\}_{i=1}^{+\infty}$ such that $r_i>\delta_{0}>0$ for all $i\in \N$, and by the definition of $\rho$, we know
\begin{align*}
	\rho\geq\int_{B_{r_i}(x_k)}e^{nv_i}\ud V_{\S^n}\geq\int_{B_{\delta_{0}}(x_k)}e^{nv_i}\ud V_{\S^n}.
\end{align*}
Taking the limit of the above inequality, we know
\begin{align*}
	\rho\geq \liminf_{i\to+\infty}\int_{B_{\delta_{0}}(x_k)}e^{nv_i}\ud V_{\S^n}\geq \frac{\alpha|\S^n|}{2}+o_{\delta_{0}}(1),
\end{align*}
if $\delta_{0}$ sufficiently small, this is impossible !

\end{pf}

We consider the exponential maps
\begin{align*}
	\exp_{x_{k,i}}: B_1(0)\subset T_{x_{k,i}}\S^n\longrightarrow \S^n,
\end{align*}
denote $\tilde{g}_i=\exp_{x_{k,i}}^{*}g_{\S^n}$ and $\tilde{v}_i(x)=v_i\circ \exp_{x_{k,i}}(x)$, then
\begin{align*}
	P_{\tilde{g}_i}\tilde{v}_i+\frac{2}{\alpha_i}Q_{g_{\S^n}}=\frac{2}{\alpha_i}Q_{g_{\S^n}}e^{n\tilde{v}_i},\qquad\qquad\mathrm{for}\qquad x\in B_1(0).
\end{align*}
The rescaling of $\tilde{v}_i$ is defined by
\begin{align}\label{hat vi}
	\hat{v}_i(x):=\tilde{v}_i(\delta_{r_i}(x))+\log r_i=\tilde{v}_i(r_ix)+\log r_i,\qquad\mathrm{for}\qquad x\in B_{1/r_i}(0).
\end{align}
where $\delta_{r_i}(x)=r_ix$. Clearly, for any $y\in B_{1/2r_i}(0)$,
\begin{align}\label{uniform volum estimate}
	\int_{B_{1/2}(y)}e^{n\hat{v}_i(x)}dx=&\int_{\delta_i(B_{1/2}(y))}e^{n\tilde{v}_i(x)}dx=(1+o_i(1))\int_{\exp_{x_{k,i}}(\delta_i(B_{1/2}(y)))}e^{nv_i}\ud V_{\S^n}\nonumber\\
	\leq&(1+o_i(1))\int_{B_{r_i}(x_{i,k})}e^{nv_i}\ud V_{\S^n}<\frac{3\alpha|\S^n|}{8}.
\end{align}
 Setting $\hat{g}_i(x)=r_i^{-2}\delta_{r_i}^{*}\tilde{g}_i=\left(\frac{2}{1+|r_ix|^2}\right)^2|dx|^2\to |dx|^2$ in $C^k_{\mathrm{loc}}(\R^n)$ sense, by the conformally invariant property of GJMS operator, then we have 
\begin{align}\label{blow up equation}
	P_{\hat{g}_i}\hat{v}_i+\frac{2r_i^n}{\alpha_i}Q_{g_{\S^n}}=\frac{2}{\alpha_i}Q_{g_{\S^n}}e^{n\hat{v}_i},\qquad\qquad\mathrm{for}\qquad x\in B_{1/r_i}(0).
\end{align}
In the following Lemma, through Lemma \ref{priori estimate lem} and \eqref{uniform volum estimate}, we arrive that if $R>1$,
\begin{align*}
	 \int_{B_{R}(0)} e^{p\hat{v}_i}\leq C(R)\qquad\mathrm{and}\qquad \int_{B_{R}(0)} |\hat{v}_i|^p\leq C(R)\qquad \mathrm{for}\qquad p=p(R)>1.
\end{align*}
Hence, the equation \eqref{blow up equation} will convergent in $C^{\infty}_{\mathrm{loc}}(\R^n)$ sense.
\begin{lem}\label{blow up lem}
	Suppose $\hat{v}_i$ is defined by \eqref{hat vi}, then
	\begin{align}
		\hat{v}_i\longrightarrow \hat{v}_{\infty}\qquad\mathrm{in}\qquad C_{\mathrm{loc}}^{\infty}(\R^n),
	\end{align}
	where 
	\begin{align*}
		\left(-\Delta\right)^{\frac{n}{2}}\hat{v}_{\infty}=\frac{2Q_{g_{\S^n}}}{\alpha}e^{n\hat{v}_{\infty}}
	\end{align*}
	and
	\begin{align*}
	\hat{v}_{\infty}(x)=\log\frac{2\lambda}{\lambda^2+|x-x_0|^2}+\frac{1}{n}\log \alpha\qquad\mathrm{and}\qquad \int_{\R^n}e^{n\hat{v}_{\infty}}=\alpha|\S^n|.
\end{align*}
\end{lem}
\begin{pf}
For any $R>0$,	we choose the cut-off function $\phi_{R}(x)$ such that $\phi_{R}(x)\equiv 1$ in $B_{R/2}(0)$ and $\phi_{R}(x)\equiv 0$ in $\R^n\backslash B_{R}(0)$, 
	\begin{align*}
		u_i=\phi_{R}	\hat{v}_i+(1-\phi_{R})a_i=a_i+\phi_{R}(\hat{v}_i-a_i) \qquad\mathrm{and}\qquad \hat{u}_i=u_i-a_i,
	\end{align*}
	where
	\begin{align*}
		a_i=\frac{1}{B_R(0)}\int_{B_R(0)}\hat{v}_i.
	\end{align*}
	We know $\hat{u}_i=\hat{v}_i-a_i$ in $B_{R/2}$ and $\hat{u}_i=0$ in $\R^n\backslash B_R(0)$. According to Lemma \ref{priori estimate lem},  for $n-kp>0$ and $k=1,2,\cdots,n-1$ we have
	\begin{align}\label{blow up: gradient estimate}
		\int_{B_{2R}(0)}|\nabla^k\hat{v}_i|^p=r_i^{kp-n}\int_{B_{2Rr_i}}|\nabla^k\tilde{v}_i|^p\leq Cr_i^{kp-n}(2Rr_i)^{n-kp}\leq CR^{n-kp}.
	\end{align}
	By the \eqref{blow up: gradient estimate} and Poincare inequality, we have
	\begin{align}\label{blow up: zero term estimate}
	\int_{B_{R}(0)}\left|\hat{v}_i-a_i\right|^p\leq C(R)\int_{B_{2R}(0)}\left|\nabla \hat{v}_i\right|^p\leq C(R).
	\end{align}

	On the other hand, $\hat{u}_i$ satisfies
	\begin{align*}
		P_{\hat{g}_i}\hat{u}_i=&	P_{\hat{g}_i}\left[\phi_{R}(\hat{v}_i-a_i)\right]=\phi_{R}P_{\hat{g}_i}\hat{v}_i+L_i(\hat{v}_i-a_i)\\
		=&\frac{2Q_{g_{\S^n}}\phi_{R}}{\alpha_i}e^{n\hat{v}_i}-\frac{2r_i^nQ_{g_{\S^n}}\phi_{R}}{\alpha_i}+L_i(\hat{v}_i-a_i):=\hat{f}_i,
	\end{align*}
	where $L_i$ is the differential operator which contains oder $0\leq k\leq n-1$ with  smooth coefficients and $\mathrm{supp}L_i\subset B_R(0)$. 
Combining \eqref{uniform volum estimate}, \eqref{blow up: gradient estimate} and \eqref{blow up: zero term estimate},  for any $y\in B_R(0)$, we deduce that
\begin{align}\label{blow up: hat fi estimate}
	\int_{B_r(y)}|\hat{f}_i|\leq& \frac{2Q_{g_{\S^n}}}{\alpha_i}\int_{B_r(y)}e^{n\hat{v}_i(x)}dx+ Cr_i^n|B_r|+ C\int_{B_r(y)}\left|L_i(\hat{v}_i-a_i)\right|\nonumber\\
	\leq&\frac{3\alpha c_n}{4\alpha_i}+Cr_i^n|B_r|+C\left(\int_{B_r(y)}\left|L_i(\hat{v}_i-a_i)\right|^p\right)^{1/p}|B_r|^{1-1/p}\nonumber\\
	\leq& \frac{3\alpha c_n}{4\alpha_i}+Cr_i^nr^n+C(R)r^{n(1-1/p)}<c_n,
\end{align}
for $r$ is sufficiently small.
   Since $\mathrm{suup} \,\,\hat{u}_i\subset B_R(0)$, we can embed the $(B_{2R}, \hat{g}_i)$ into torus $(\T^n,h_i)$ such that $h_i=\hat{g}_i$ in $B_{2R}(0)$  in $T^n\backslash B_{3R}$ and $h_i\to |dx|^2$ as $i\to+\infty$. Through Lemma \ref{concentration compactness lem} and \eqref{blow up: hat fi estimate}, we know there exist $p:=p(R)>n$ such that 
   \begin{align}\label{blow up:exp hat ui estimate}
   	\int_{B_{R}(0)} e^{p\hat{u}_i}\leq C(R).
   \end{align}

   We claim that:
   \begin{align}\label{blow up: a i estimate}
   	|a_i|\leq C(R).
   \end{align}
   On the one hand, Jensen inequality implies
   \begin{align}\label{blow up: ai upper bound}
   	e^{na_i}\leq& \int_{B_{R}(0)}e^{n\hat{v}_i}=\int_{B_{Rr_i}(0)}e^{n\tilde{v}_i}\nonumber\\
   	=&(1+o_i(1))\int_{\exp_{x_{k,i}}(B_{Rr_i}(0))}e^{nv_i}\ud V_{\S^n}\leq C.
   \end{align}
   On the other hand, 
   \begin{align}\label{blow up: lower bound}
   	\int_{B_{R/2}(0)}e^{n\hat{v}_i}=&(1+o_i(1))\int_{\exp_{x_{k,i}}(B_{Rr_i/2}(0))}e^{nv_i}\ud V_{\S^n}\nonumber\\
   	\geq&(1+o_i(1))\int_{B_{r_i}(x_{k,i})}e^{nv_i}\ud V_{\S^n}\geq\frac{\rho}{2}
   \end{align}
    and 
    \begin{align}\label{blow up: up bound}
    	\int_{B_{R/2}(0)}e^{n\hat{v}_i}=&e^{na_i}\int_{B_{R/2}(0)}e^{n\hat{u}_i}\leq e^{na_i}	\left(\int_{B_{R/2}(0)} e^{p\hat{u}_i}\right)^{1/p}|B_{R/2}|^{1-1/p}\nonumber\\
    	\leq& C(R)e^{na_i},
    \end{align}
    where \eqref{blow up: up bound} follows by \eqref{blow up:exp hat ui estimate}. Then,
    \eqref{blow up: lower bound} and \eqref{blow up: up bound} will imply 
    \begin{align} \label{blow up: ai lower  bound} 
    	a_i\geq -C(R).
    \end{align}
    We complete the proof the claim  by \eqref{blow up: ai upper bound} and \eqref{blow up: ai lower  bound} .
     Hence, \eqref{blow up:exp hat ui estimate} and \eqref{blow up: a i estimate} will yield
     \begin{align}\label{blow up:exp hat vi upper bound}
     \int_{B_{R/2}(0)} e^{p\hat{v}_i}=e^{pa_i}	\int_{B_{R/2}(0)} e^{p\hat{u}_i}\leq C(R)
     \end{align}
      and combining with \eqref{blow up: zero term estimate} we obtain
      \begin{align}\label{blow up: hat vi upper bound}
      	 \int_{B_{R/2}(0)} |\hat{v}_i|^p\leq C(R).
      \end{align}
     Recall the equation
     \begin{align*}
     	P_{\hat{g}_i}\hat{v}_i+\frac{2r_i^n}{\alpha_i}Q_{g_{\S^n}}=\frac{2}{\alpha_i}Q_{g_{\S^n}}e^{n\hat{v}_i},\qquad\qquad\mathrm{for}\qquad x\in B_{1/r_i}(0),
     \end{align*}
     we notice that \eqref{blow up:exp hat vi upper bound} , \eqref{blow up: hat vi upper bound} and high order elliptic estimate (see the Lemma 7.6 in \cite{Chen Lu Zhu}) give that
     \begin{align*}
||\hat{v}_i||_{C^{\beta}(B_{R/4}(0))}\leq C||\hat{v}_i||_{W^{n,p}(B_{R/4}(0))}\leq C(R),
     \end{align*}
     for some $\beta:=\beta(R)\in (0,1)$.
     Again, elliptic estimate also gives
     \begin{align}
     	||\hat{v}_i||_{C^{n,\beta}(B_{R/4}(0))}\leq C(R).
     \end{align}
     Finally, we establish
     \begin{align*}
     	||\hat{v}_i||_{C^{m}(B_{R/4}(0))}\leq C(m,R)\qquad\mathrm{for\,\,\,any}\quad m\in\N,\quad R>0.
     \end{align*}
     So, we have \begin{align*}
     \hat{v}_i\longrightarrow \hat{v}_{\infty}\qquad\mathrm{in}\qquad C_{\mathrm{loc}}^{\infty}(\R^n)
     \end{align*}
     and
     where 
     \begin{align*}
     \left(-\Delta\right)^{\frac{n}{2}}\hat{v}_{\infty}=\frac{2Q_{g_{\S^n}}}{\alpha}e^{n\hat{v}_{\infty}}.
     \end{align*}
     Noticing that \eqref{blow up: gradient estimate}, we have 
     \begin{align}\label{blow up: }
     	\int_{B_{R}(0)}|\nabla^k\hat{v}_{\infty}|\leq CR^{n-k}\qquad\mathrm{for}\qquad k=1,2,\cdots,n-1.
     \end{align}
     Theorem \ref{classification thm } tell us that
     \begin{align*}
     	\hat{v}_{\infty}(x)=\log\frac{2\lambda}{\lambda^2+|x-x_0|^2}+\frac{1}{n}\log \alpha\qquad\mathrm{and}\qquad \int_{\R^n}e^{n\hat{v}_{\infty}}=\alpha|\S^n|.
     \end{align*}
\end{pf}

\begin{cor}\label{cor lower bound}
	With the same notation as Lemma \ref{blow up lem}, up to a subsequence, for any $r>0$ there holds
	\begin{align*}
		\liminf_{i\to+\infty}\int_{B_{r}(x_k^{'})}e^{nv_i}\geq \alpha|\S^n|.
	\end{align*}
\end{cor}
\begin{pf}
	This corollary follows by Lemma \ref{blow up lem} and for any $R>1$
	\begin{align*}
		\int_{B_{R}(0)}e^{n\hat{v}_i}=&(1+o_i(1))\int_{\exp_{x_{k,i}}(B_{Rr_i}(0))}e^{nv_i}\ud V_{\S^n}\leq(1+o_i(1))\int_{B_{r}(x_k^{'})}e^{nv_i}\ud V_{\S^n}.
	\end{align*}
\end{pf}

\begin{lem}\label{blow up number lem}
	Suppose $\{x_k\}_{k=1}^{j_0}$ are the all concentration points of $\{v_i\}_{i=1}^{\infty}$, then
	\begin{align*}
		x_k=x_k^{'}\qquad\mathrm{for}\qquad k=1,2,\cdots, j_0
	\end{align*}
	and $j_0=1$.
\end{lem}
\begin{pf}
	Suppose there exist $k$ such that $x_k\not=x_k^{'}$, then for small $r>0$, up to a subsequence, 
	\begin{align*}
		\liminf_{i\to+\infty}\int_{B_{r}(x_k^{'})}|f_i|=&\liminf_{i\to+\infty}\frac{2Q_{g_{\S^n}}}{\alpha_i}\int_{B_{r}(x_k^{'})}e^{nv_i}+o_{r}(1)\\
		\geq&\liminf_{i\to+\infty}\frac{2Q_{g_{\S^n}}}{\alpha_i}\int_{B_{r_i}(x_{k,i})}e^{nv_i}\ud V_{\S^n}+o_{r}(1)\\
		\geq&2Q_{g_{\S^n}}|\S^n|+o_r(1)>c_n.
	\end{align*}
	Then, $x_k^{'}$ is also a concentration point, $\{v_i\}_{i=1}^{\infty}$ has $j_0+1$ concentration points, this is a contradiction. For any $1\leq k\leq j_0$, then
	\begin{align*}
		\liminf_{i\to+\infty}\int_{B_{r}(x_k)}e^{nv_i}\ud V_{\S^n}\geq \alpha|\S^n|.
	\end{align*}
	Hence, 
	\begin{align*}
		|\S^n|=\liminf_{i\to+\infty}\sum_{k=1}^{j_0}\int_{B_{r}(x_k)}e^{nv_i}\ud V_{\S^n}\geq j_0\alpha|\S^n|,
	\end{align*}
	$\alpha\in (1/2, 1]$ will imply $j_0=1$.
\end{pf}

\begin{lem}\label{one point volume vanishing lem}
		Suppose $\{v_i\}_{i=1}^{+\infty}$ has only one concentration point $x_1$, up to a subsequence,  then 
		for any $r>0$
		\begin{align}\label{concentration point concentrate}
				\liminf_{i\to+\infty}\int_{B_{r}(x_1)}e^{nv_i}\ud V_{\S^n}\geq \alpha|\S^n|
		\end{align}
		and
		\begin{align}\label{concentration point out vanish}
		\lim_{i\to+\infty}\int_{\S^n\backslash B_{r}(x_1)}e^{nv_i}\ud V_{\S^n}=0.
		\end{align} 
\end{lem}
\begin{pf}
	For \eqref{concentration point concentrate}, it is easy to see by Corollary \ref{cor lower bound}  and Lemma \ref{blow up number lem}. 
	We assume that \eqref{concentration point out vanish} is not true, there exist $r_0$ such that
	\begin{align}\label{one point volume vanishing a1}
		\int_{\S^n\backslash B_{r_0}(x_1)}e^{nv_i}\ud V_{\S^n}\geq c_0>0.
	\end{align}
	Since $\{v_i\}_{i=1}^{+\infty}$ has only one concentration point $x_1$, then for any $y\in \S^n\backslash B_{r_0}(x_1)$, there exist $\delta_{y}>0$ and $r_{y}>0$ which satisfy
	\begin{align*}
		\int_{B_{r_{y}}(y)}|f_i|\ud V_{\S^n}< c_n-\delta_{y}, \qquad\mathrm{for}\qquad i\gg 1.
	\end{align*}
	By the same proof as for the Lemma \ref{concentration compactness lem}, you can obtain that there exist $p>n$ and a subsequence $\{v_i\}_{i=1}^{+\infty}$ such that
	\begin{align}\label{one point volume vanishing a2}
		\int_{\S^n\backslash B_{r_0}(x_1)}e^{p\left(v_i-\bar{v}_i\right)}\ud V_{\S^n}\leq C.
	\end{align}
	Then \eqref{one point volume vanishing a1} and \eqref{one point volume vanishing a2} will yield
	\begin{align}\label{one point volume vanishing a3}
		c_0\leq \int_{\S^n\backslash B_{r_0}(x_1)}e^{nv_i}\ud V_{\S^n}=e^{n\bar{v}_i}\int_{\S^n\backslash B_{r_0}(x_1)}e^{n\left(v_i-\bar{v}_i\right)}\ud V_{\S^n}\leq Ce^{n\bar{v}_i}.
	\end{align}
	In the next step, we claim that
	\begin{align}\label{one point volume vanishing a4}
		\bar{v}_i\longrightarrow-\infty\qquad\mathrm{as}\qquad i\longrightarrow+\infty.
	\end{align}
	For any $r>0$
	\begin{align*}
		\liminf_{i\to+\infty}\int_{B_{r}(x_1)}e^{nv_i}\ud V_{\S^n}\geq \alpha|\S^n|,
	\end{align*}
	and 
	\begin{align*}
		\left|G_{\S^n}(x,y)-\frac{1}{c_n}\log\frac{1}{|x-y|}\right|\leq C\qquad\mathrm{for}\qquad x\not=y\in\S^n.
	\end{align*}
	Therefore, 
	\begin{align}\label{one point volume vanishing a5}
		v_i(x)-\bar{v}_i=&\frac{2Q_{g_{\S^n}}}{\alpha_i}\int_{\S^n}G_{\S^n}(x,y)\left(e^{nv_i(y)}-1\right)\ud V_{\S^n}(y)\nonumber\\
		\geq& \frac{2Q_{g_{\S^n}}}{\alpha_i}\left(\int_{B_r(x_1)}+\int_{\S^n\backslash B_{r}(x_1)}G_{\S^n}(x,y)e^{nv_i(y)}\ud V_{\S^n}(y)\right)-C(n)\nonumber\\
		\geq&\frac{2Q_{g_{\S^n}}}{\alpha_i}\int_{ B_{r}(x_1)}G_{\S^n}(x,y)e^{nv_i(y)}\ud V_{\S^n}(y)-C(n),
	\end{align}
	the last inequality follows by 
	\begin{align*}
	G_{\S^n}(x,y)\geq \frac{1}{c_n}\log\frac{1}{|x-y|}-C\geq  \frac{1}{c_n}\log\frac{1}{\pi}-C.
	\end{align*}
	Suppose $x\in \S^n\backslash B_{2r}(x_1)$, then 
	\begin{align}\label{one point volume vanishing a6}
		\left|\log\frac{1}{|x-y|}-\log\frac{1}{|x-x_1|}\right|\leq C(n).
	\end{align}
	Hence, from \eqref{one point volume vanishing a5} and \eqref{one point volume vanishing a6} we know
	\begin{align*}
		e^{nv_i(x)}\geq& Ce^{n\bar{v}_i}\exp\left(\frac{2nQ_{g_{\S^n}}}{c_n\alpha_i}\log\frac{1}{|x-x_1|}\int_{B_{r}(x_1)}e^{nv_i}\ud V_{\S^n}\right)\\
		\geq& Ce^{n\bar{v}_i}\frac{1}{|x-x_1|^{2n-\e}}\qquad\mathrm{for}\qquad i>N(r)\gg 1\qquad \mathrm{and}\qquad x\in \S^n\backslash B_{2r}(x_1).
	\end{align*}
	However, this will cause
\begin{align*}
	|\S^n|=&\int_{\S^n}e^{nv_i}\ud V_{\S^n}\geq \int_{B_{1}(x_1)\backslash B_{2r}(x_1) }e^{nv_i}\ud V_{\S^n}\\
	\geq& Ce^{n\bar{v}_i}\int_{B_{1}(x_1)\backslash B_{2r}(x_1) }\frac{1}{|x-x_1|^{2n-\e}}\ud V_{\S^n}(x)\\
	=& Ce^{n\bar{v}_i}\frac{1}{r^{n-\e}}\qquad\mathrm{for}\qquad i\gg 1.
\end{align*}
Due the arbitrariness  of $r$, we know $\bar{v}_i\to -\infty$ as $i\to+\infty$. The final result is given by \eqref{one point volume vanishing a3} and \eqref{one point volume vanishing a4}.
\end{pf}

\begin{lem}\label{sec compactness thm}
	Suppose $\{v_i\}_{i=1}^{+\infty}$ has only one concentration point $x_1$ and satisfies
	\begin{align*}
	P_{\S^n}v_{i}+\frac{2}{\alpha}Q_{g_{\S^n}}=\frac{2}{\alpha}Q_{g_{\S^n}}e^{nv_{i}},
	\end{align*}
	then $\{v_i\}_{i=1}^{+\infty}$ are uniformly bounded in $C^{\infty}(\S^n)$ sense.
\end{lem}
\begin{pf}
	If not, then the second case in Lemma \ref{concentration compactness lem} will occur.
	Firstly, we claim that: up to a subsequence,
	for any $\phi\in C^{\infty}(\S^n)$ and $\phi(x_1)>0$, 
	\begin{align}\label{sec compactness thm: claim}
	\liminf_{i\to+\infty}\fint_{\S^n}e^{nv_i}\phi \ud V_{\S^n}\geq\alpha\phi(x_1).
	\end{align}
	We note that for sufficiently small  $r>0$,
	\begin{align}\label{sec compactness thm: a1}
		&\int_{\S^n}e^{nv_i}\phi \ud V_{\S^n}\nonumber\\		
		=&\int_{ B_{r}(x_1)}e^{nv_i}\phi \ud V_{\S^n}+\int_{\S^n\backslash B_{r}(x_1)}e^{nv_i}\phi \ud V_{\S^n}\nonumber\\
		=&\phi(x_1)\int_{ B_{r}(x_1)}e^{nv_i}\ud V_{\S^n}+\int_{ B_{r}(x_1)}e^{nv_i}\left(\phi(x)-\phi(x_1)\right) \ud V_{\S^n}+\int_{\S^n\backslash B_{r}(x_1)}e^{nv_i}\phi \ud V_{\S^n}\nonumber\\
		=&\phi(x_1)\int_{ B_{r}(x_1)}e^{nv_i}\ud V_{\S^n}+I+II,
	\end{align}
	then
		\begin{align}\label{sec compactness thm: a2}
	|I|=&\left|\fint_{\S^n}e^{nv_i(x)}\left(\phi(x)-\phi(x_1)\right) \ud V_{\S^n}(x)\right|\nonumber\\
	\leq&C||\nabla\phi||_{L^{\infty}(\S^n)}r\fint_{B_r(x_1)}e^{nv_i(x)}\ud V_{\S^n}(x)\leq Cr.
	\end{align}
	The Lemma \ref{one point volume vanishing lem} implies that
	\begin{align}\label{sec compactness thm: a3}
	|II|\leq C(n)||\phi||_{L^{\infty}(\S^n)}\int_{\S^n\backslash B_{r}(x_1)}e^{nv_i} \ud V_{\S^n}\to 0\qquad\mathrm{as}\qquad i\to+\infty,
	\end{align}
	 from \eqref{sec compactness thm: a1},  \eqref{sec compactness thm: a2} and \eqref{sec compactness thm: a3} we obtain
	\begin{align*}
		\liminf_{i\to+\infty}\fint_{\S^n}e^{nv_i}\phi \ud V_{\S^n}\geq& \liminf_{i\to+\infty}\phi(x_1)\int_{ B_{r}(x_1)}e^{nv_i}\ud V_{\S^n}-Cr\\
		\geq& \alpha|\S^n|\phi(x_1)-Cr.
	\end{align*}
	Letting $r\to 0$, we get the \eqref{sec compactness thm: claim}. We can suppose $x_1$ as the north pole for $\S^n$, letting $\phi=x^{n+1}$,  \eqref{sec compactness thm: claim} will give 
	\begin{align*}
	\liminf_{i\to+\infty}\fint_{\S^n}e^{nv_i}x^{n+1} \ud V_{\S^n}\geq \alpha.
	\end{align*}
	This is a contradiction with Kazdan-Warner condition \eqref{direc assu} in Lemma \ref{Kazdan-Warner lem}.
\end{pf}

\begin{thm}
Suppose $n=2m\geq 4$, then for any $\e>0$ the following set
\begin{align*}
\mathcal{S}_{\e}=\left\{v_{\alpha}\Big|\alpha P_{\S^n}v_{\alpha}+2Q_{g_{\S^n}}=2Q_{g_{\S^n}}e^{nv_{\alpha}}\quad\mathrm{and}\quad \alpha\in \left(\frac{1}{2}+\e,1\right)\right\}
\end{align*}
is compact in $C^{\infty}(\S^n)$ topology.
\end{thm}
\begin{pf}
	If not, there exist $\alpha_i\to\alpha$ such that $||v_{\alpha_i}||_{C^0(\S^n)}\to+\infty$, then any subsequence of $\{v_{\alpha_i}\}$ will also blow up. According to  Lemma \ref{concentration compactness lem} and     Lemma \ref{compactness lem}, we know the second case in Lemma \ref{concentration compactness lem}  occurs, i.e, there exist finite concentration points (Here we can choose a subsequence of $\{v_{\alpha_i}\}$, which we still call it by $\{v_{\alpha_i}\}$, such that the number of concentration point is maximal) for $\{v_{\alpha_i}\}$. But uo to a subsequence,  the Lemma \ref{blow up number lem} tells us that there exists at most one concentration point, and finally Lemma \ref{sec compactness thm} will imply $||v_{\alpha_i}||_{C^0(\S^n)}\leq C$ for  $i\gg1$. This is contradictory to our assumption.
	
\end{pf}

	\noindent \textbf{Proof of Theorem \ref{main thm 2.1}:}

	If not, there exist $\alpha_i\to 1$ such that $v_{\alpha_i}\not\equiv 0$. We also note that  
\begin{align*}
v_{\alpha_i}\not\equiv \mathrm{const}\Longleftrightarrow v_{\alpha_i}\not\equiv 0.
\end{align*}
By the Theorem \ref{main thm 2}, we know $v_{\alpha_i}\to v$ in $C^{\infty}(\S^n)$ and $v$ satisfies
\begin{align}\label{uniqueness thm formula a}
P_{\S^n}v+2Q_{g_{\S^n}}=2Q_{g_{\S^n}}e^{nv}
\end{align}
and
\begin{align}\label{uniqueness thm formula a1}
\int_{\S^n}x^lv\ud V_{\S^n}=0\quad\mathrm{and}\quad	\int_{\S^n}x^le^{nv}\ud V_{\S^n}=0,\quad\mathrm{for}\quad l=1,2,\cdots,n+1.
\end{align}
We claim that
\begin{align*}
v\equiv 0.
\end{align*}
For this purpose, we choose a critical point $p$ of $v$, let $u(x)=v\circ I^{'}(x)+\log\frac{2}{1+|x|^2}$, where $I^{'}:\R^n \to \S^n\backslash\{-p\}$, then we know $0$ is a critical point of $u(x)$. We obtain that $u$ satisfies
\begin{align*}
\left(-\Delta\right)^{\frac{n}{2}}u=(n-1)!e^{nu}\qquad\mathrm{in}\qquad\R^n,
\end{align*}
by the classification of high order conformally invariant equation (see \cite{Chang and Yang} or \cite{Wei and Xu 1}), we get
\begin{align*}
u(x)=\log\frac{2\lambda}{1+\lambda^2|x|^2}.
\end{align*}
Clearly, there holds
\begin{align*}
\int_{\S^n}x^{n+1}e^{nv}\ud V_{\S^n}=&\int_{\R^n}e^{nu(x)}\frac{|x|^2-1}{|x|^2+1}dx=\int_{\R^n}\left(\frac{2\lambda}{1+\lambda^2|x|^2}\right)^n\frac{|x|^2-1}{|x|^2+1}dx\\
=&|\S^{n-1}|\int_{0}^{+\infty}\left(\frac{2\lambda }{1+\lambda^2r^2}\right)^n\frac{r^2-1}{r^2+1}r^{n-1}dr\xrightarrow{s=\lambda r}\\
=&|\S^{n-1}|\int_{0}^{+\infty}\left(\frac{2}{1+s^2}\right)^ns^{n-1}\frac{s^2-\lambda^2}{s^2+\lambda^2}ds:=C(\lambda).
\end{align*}
You can check that $C(\lambda)$ is strictly decreasing above $\lambda$. Hence, \eqref{uniqueness thm formula a1} imply that $C(\lambda)=0$, but we know $C(1)=0$, then $\lambda=1$ and $v\equiv 0$.
Letting $\bar{v}_{\alpha_i}=\fint_{\S^n}v_{\alpha_i}\ud V_{\S^n}$ and $\mathring{v}_i=v_{\alpha_i}-\bar{v}_{\alpha_i}\not\equiv 0$. Here, we adopt the strategy that first appeared on page 228 of Chang and Yang's seminal work \cite{Chang and Yang acta}.  Since
\begin{align*}
\int_{\S^n}x^l\mathring{v}_i\ud V_{\S^n}=0\qquad\mathrm{and}\qquad\int_{\S^n}\mathring{v}_i\ud V_{\S^n}=0,
\end{align*}
it mans that $\mathring{v}_i\bot\mathcal{H}_1$, where $\mathcal{H}_1=\mathrm{span}\{1, x^{1},\cdots,x^{n+1}\}$. Since $\lambda_2(P_{\S^n})=(n+1)!$, applying the min-max principle yields 
\begin{align}\label{uniqueness thm formula f}
(n+1)!\fint_{\S^n}\mathring{v}_i^2\ud V_{\S^n}\leq& \fint_{\S^n}\mathring{v}_i P_{\S^n}\mathring{v}_i\ud V_{\S^n}.
\end{align}
On the other hand, we have 
\begin{align*}
P_{\S^n}\mathring{v}_i=P_{\S^n}v_{\alpha_i}-\bar{v}_{\alpha_i}P_{\S^n}(1)=\frac{2Q_{\S^n}}{\alpha_i}\left(e^{nv_{\alpha_i}}-1\right)-\bar{v}_{\alpha_i}P_{\S^n}(1).
\end{align*}
Thus, it follows that  
\begin{align}\label{uniqueness thm formula g}
\fint_{\S^n}\mathring{v}_i P_{\S^n}\mathring{v}_i\ud V_{\S^n}=&\frac{2Q_{\S^n}}{\alpha_i}\int_{\S^n}\mathring{v}_ie^{nv_{\alpha_i}}\ud V_{\S^n}\nonumber\\
=&\frac{2Q_{g_{\S^n}}}{\alpha_i}e^{n\bar{v}_{\alpha_i}}\int_{\S^n}\mathring{v}_ie^{n\mathring{v}_i}\ud V_{\S^n}\nonumber\\
=&\frac{2Q_{g_{\S^n}}}{\alpha_i}e^{n\bar{v}_{\alpha_i}}\fint_{\S^n}\left(e^{n\mathring{v}_i}-1\right)\mathring{v}_i\ud V_{\S^n}\nonumber\\
=&n!(1+o_i(1))\fint_{\S^n}\mathring{v}_i^2\ud V_{\S^n}
\end{align}
the last equality follows by  $\bar{v}_{\alpha_i}\to 0$ , $\mathring{v}_i\to 0$ uniformly and
\begin{align*}
e^{n\mathring{v}_i}-1=(1+o_i(1))n\mathring{v}_i.
\end{align*} 
From \eqref{uniqueness thm formula f} and \eqref{uniqueness thm formula g}, we obtain $\mathring{v}_i\equiv 0$ for $i\gg1$, this is a contradiction !\\

	\noindent \textbf{Proof of Theorem \ref{main thm 2.2}:} 	
	By repeating the above proof and replacing $\alpha_i \to 1$ with $\alpha_i \to \alpha_0 \in \Gamma$, the result can be obtained.

\subsection{critical case}\label{sec 4.4}

At last, we consider the critical case, namely, 
\begin{align}
		P_{\S^n}v_{\alpha_i}+\frac{2}{\alpha_i}Q_{g_{\S^n}}=\frac{2}{\alpha_i}Q_{g_{\S^n}}e^{nv_{\alpha_i}}\qquad\mathrm{as}\qquad \alpha_i\to\frac{1}{2}.
\end{align}
Through Lemma \ref{sec compactness thm}, we know if $\{v_i\}_{i=1}^{+\infty}$ blow up, then $\{v_i\}_{i=1}^{+\infty}$ has at least two concentration points $\{x_1, x_2\cdots,x_{j_0}\}\subset\S^n$. But according to the Corollary \ref{cor lower bound}, we know
	\begin{align}\label{critical volum}
		\liminf_{i\to+\infty}\int_{B_r(x_k)}e^{nv_i}\ud V_{\S^n}\geq\frac{1}{2}|\S^n|\qquad\mathrm{for\,\,\,\,\,any}\qquad r>0.
	\end{align}
Since $\int_{\S^n}e^{nv_i}\ud V_{\S^n}=|\S^n|$, then $j_0\leq2$. In this case, $\{v\}_{i=1}^{+\infty}$ blow up iff $\{v\}_{i=1}^{+\infty}$ has  two concentration points $\{x_1, x_2\}\subset\S^n$.  We will prove for any $\phi\in C^{\infty}(\S^n)$,
\begin{align*}
	\lim_{i\to+\infty}\fint_{\S^n}e^{nv_i}\phi\ud V_{\S^n}=\frac{\phi(x_1)+\phi(x_2)}{2}.
\end{align*}
	Combining with Kazdan-Warner condition, we know 
	
		\begin{align*}
		\int_{\S^n}x^{n+1}e^{nv_{i}}\ud V_{\S^n}=0,
		\end{align*}
	we obtain $\{x_1, x_2\}=\{x_1, -x_1\}\subset\S^n$.

	\begin{thm}\label{critical blow up thm}
		Suppose $v_i=v_{\alpha_i}$ is the solution of 
		\begin{align}
		P_{\S^n}v_{\alpha_i}+\frac{2}{\alpha_i}Q_{g_{\S^n}}=\frac{2}{\alpha_i}Q_{g_{\S^n}}e^{nv_{\alpha_i}}\qquad\mathrm{as}\qquad \alpha_i\to\frac{1}{2}.
		\end{align}
		Then $\{v_i\}_{i=1}^{+\infty}$ either is uniformly bounded or up to a subsequence, it has two concentration points $\{x_1, -x_1\}\subset\S^n$ and
	\begin{align*}
		e^{nv_i}\ud V_{\S^n}\rightharpoonup \frac{\delta_{x_1}+\delta_{-x_1}}{2}|\S^n|.
	\end{align*}
	\end{thm}
	\begin{pf}
		By the Lemma \ref{compactness lem} and Lemma  \ref{sec compactness thm}, if  $\{v\}_{i=1}^{+\infty}$ blow up, then $\{v\}_{i=1}^{+\infty}$ has at least two concentration points $\{x_1, x_2\cdots,x_{j_0}\}\subset\S^n$. For any $1\leq k\leq j_0$, the Lemma \ref{blow up lem} implies that, up to a subsequence, 
		\begin{align}\label{critical blow up thm a_1}
		\liminf_{i\to+\infty}\int_{B_r(x_k)}e^{nv_i}\ud V_{\S^n}\geq\frac{1}{2}|\S^n|\qquad\mathrm{for\,\,\,\,\,any}\qquad r>0.
		\end{align}
		And,
		\begin{align}\label{critical blow up thm a_2}
		\limsup_{i\to+\infty}	\sum_{k=1}^{j_0}\int_{B_r(x_k)}e^{nv_i}\ud V_{\S^n}\leq 	\limsup_{i\to+\infty}\int_{\S^n}e^{nv_i}\ud V_{\S^n}=|\S^n|.
		\end{align}
		By \eqref{critical blow up thm a_1} and \eqref{critical blow up thm a_2}, we know $j_0=2$, 
		\begin{align}\label{critical blow up thm a_3}
			\lim_{i\to+\infty}\int_{B_r(x_1)}e^{nv_i}\ud V_{\S^n}=	\lim_{i\to+\infty}\int_{B_r(x_2)}e^{nv_i}\ud V_{\S^n}=\frac{1}{2}|\S^n|\qquad\mathrm{for\,\,\,\,\,any}\qquad r>0
		\end{align}
		and 
		\begin{align}\label{critical blow up thm a_4}
				\lim_{i\to+\infty}\int_{\S^n\backslash \left(B_r(x_1)\cup  B_r(x_2)\right)}e^{nv_i}\ud V_{\S^n}=0
				\qquad\mathrm{for\,\,\,\,\,any}\qquad r>0.
		\end{align}
		We claim that: for any $\phi\in C^{\infty}(\S^n)$,
		\begin{align}\label{critical blow up thm claim}
		\lim_{i\to+\infty}\fint_{\S^n}e^{nv_i}\phi\ud V_{\S^n}=\frac{\phi(x_1)+\phi(x_2)}{2}.
		\end{align}
		Since
		\begin{align*}
			&\left|\fint_{\S^n}e^{nv_i}\phi\ud V_{\S^n}-\frac{\phi(x_1)+\phi(x_2)}{2}\right|\\
			\leq&\left|\fint_{B_r(x_1)}e^{nv_i}\left(\phi(x)-\phi(x_1)\right)\ud V_{\S^n}\right|+\left|\fint_{B_r(x_2)}e^{nv_i}\left(\phi(x)-\phi(x_2)\right)\ud V_{\S^n}\right|\\
			+&||\phi||_{L^{\infty}(\S^n)}\left|\fint_{B_r(x_1)}e^{nv_i}\ud V_{\S^n}-\frac{1}{2}\right|+||\phi||_{L^{\infty}(\S^n)}\left|\fint_{B_r(x_2)}e^{nv_i}\ud V_{\S^n}-\frac{1}{2}\right|\\
			+&||\phi||_{L^{\infty}(\S^n)}\int_{\S^n\backslash B_r(x_1)\cup  B_r(x_2)}e^{nv_i}\ud V_{\S^n}\\
			\leq& C||\nabla\phi||_{L^{\infty}(\S^n)}r+o_i(1)\qquad\mathrm{for\,\,\,\,\,any}\qquad r>0.
		\end{align*}
		the last inequality follows by \eqref{critical blow up thm a_3} and \eqref{critical blow up thm a_4}. Taking limit for $i$ and letting $r\to 0$, we get the \eqref{critical blow up thm claim}.
		We can suppose $x_1$ as the north pole for $\S^n$, letting $\phi=x^{n+1}$, \ref{critical blow up thm claim} and Kazdan-Warner condition will give  
		\begin{align*}
		\lim_{i\to+\infty}\fint_{\S^n}e^{nv_i}x^{n+1} \ud V_{\S^n}=\frac{1+x^{n+1}(x_2)}{2}=0,
		\end{align*}
		hence $\{x_1, x_2\}=\{x_1, -x_1\}\subset\S^n$.
	\end{pf}

Finally, we point out that $v_i\to -\infty$ in $\S^n\backslash \{x_1, -x_1\}$, this phenomenon  is very common, which also occurs in two dimensions see \cite{Brezis,Liyanyan}.

	\begin{lem}
		With the same assumption as Theorem \ref{critical blow up thm}, $\{v_i\}$ blow up at $\{x_1, -x_1\}\subset\S^n$, then
		\begin{align}\label{vi limit formula a}
		v_i(x)-\bar{v}_i\longrightarrow 4Q_{g_{\S^n}}|\S^n|\left(\frac{G_{\S^n}(x, x_1)+G_{\S^n}(x, -x_1)}{2}-\fint_{\S^n}G_{\S^n}(x,y   \ud V_{{\S^n}}(y)\right)
		\end{align}
		in $ C^{\infty}\left(\S^n\backslash \{x_1, -x_1\}\right)$.	Especially, $v_i\to-\infty$ in $\S^n\backslash \{x_1, -x_1\}$.
	\end{lem}
	\begin{pf}
By the	Theorem \ref{critical blow up thm} and 
		\begin{align*}
				v_i(x)-\bar{v}_i=\frac{2Q_{g_{\S^n}}}{\alpha_i}\int_{\S^n}G_{\S^n}(x,y)\left(e^{nv_i(y)}-1\right)\ud V_{\S^n}(y),
		\end{align*}
	letting $i\to+\infty$, we have
	\begin{align*}
		v_i(x)-\bar{v}_i\longrightarrow 2Q_{g_{\S^n}}|\S^n|\left(G_{\S^n}(x, x_1)+G_{\S^n}(x, -x_1)\right)-d_n\qquad\mathrm{in}\qquad C^{\infty}\left(\S^n\backslash \{x_1, -x_1\}\right).
	\end{align*}
	 We also point out that $d_n$ is independent of $x$. For any $x^{'}$, there exist $A\in SO(n)$ such that $x^{'}=Ax$, then
	\begin{align*}
	d_n(x^{'})=&4Q_{g_{\S^n}}\int_{\S^n}G_{\S^n}(Ax,y)\ud V_{\S^n}(y)=4Q_{g_{\S^n}}\int_{\S^n}G_{\S^n}(x,Ay)\ud V_{\S^n}(y)\\
	=&4Q_{g_{\S^n}}\int_{\S^n}G_{\S^n}(x,z)\ud V_{\S^n}(z)=d_n(x),
	\end{align*}
	the third equality follows by $\left(A^{-1}\right)^{*}\ud V_{\S^n}=\ud V_{\S^n}$. By \eqref{one point volume vanishing a4}, we know $\bar{v}_i\to-\infty$, so we obtain $v_i(x)\to-\infty$ in $\S^n\backslash \{x_1, -x_1\}$.
	\end{pf}

	{\noindent\small{\bf Acknowledgment:} We deeply appreciate the referee's invaluable feedback and insightful suggestions, which have greatly enhanced the quality of this paper.}

		{\noindent\small{\bf Data availability:} Data sharing not applicable to this article as no datasets were generated or analysed during the current study.
		\section*{Declarations}
		{\noindent\small{\bf Conflict of interest:} On behalf of all authors, the corresponding author states that there is no conflict of interest.

		\bibliographystyle{unsrt}

		\bigskip
		
		\noindent S. Zhang
		
		\noindent School of Mathematics, Nanjing University, \\
		Nanjing 210093, China\\[1mm]
		Email: \textsf{dg21210019@smail.nju.edu.cn}
		
		\medskip  	
		
	\end{document}